\documentclass[11pt]{article}
\usepackage{amsfonts}
\usepackage{lscape}
\usepackage{latexsym,amsmath}
\usepackage{amssymb,array}
\usepackage{xcolor}
\makeatletter \oddsidemargin 0in \evensidemargin 0in \textwidth
16cm \RequirePackage[dvips]{graphicx} \textheight 20cm
\newtheorem{t1}{Theorem}[section]
\newtheorem{p1}{Proposition}[section]

\newtheorem{d1}{Definition}[section]
\newtheorem{r1}{Remark}[section]

\begin{document}
	\title{Fractional order entropy-based decision-making models under risk}
	\author{Poulami Paul\and Chanchal Kundu\footnote{Corresponding
			author e-mail: chanchal$_{-}$kundu@yahoo.com/ckundu@rgipt.ac.in}
		\and\and
		Department of Mathematical Sciences\\
		Rajiv Gandhi Institute of Petroleum Technology\\
		Jais 229 304, U.P., India}
	\date{Revised version to appear in \textit{Communications in Nonlinear Science and}\\ \textit{Numerical Simulation}, \textcopyright{} by Elsevier. \\ Submitted: November, 2024} 
\maketitle
\begin{abstract} 
The construction of an efficient portfolio with a good level of return and minimal risk depends on selecting the optimal combination of stocks. This paper introduces a novel decision-making framework for stock selection based on fractional-order entropy due to Ubriaco. By tuning the fractional parameter $ \alpha $, the model captures varying attitudes of individuals toward risk. Values of $ \alpha $ near 1 indicate high risk tolerance (adventurous attitude), while those near 0 reflect risk aversion (conservative attitude). The sensitivity of the fractional order entropy to changing risk preferences of decision-makers is demonstrated through four real-world portfolio models, namely, large-cap, mid-cap, diversified, and hypothetical. Furthermore, two new risk measures, termed as expected utility–fractional entropy (EU-FE) and expected utility–fractional entropy and variance (EU-FEV), are introduced to develop decision models aligned with investors' risk preferences. The effectiveness of the decision model is further tested with financial stock market data of PSI 20 index by finding efficient frontiers of portfolio with the aid of artificial neural network.
\end{abstract}
\section{Introduction}
The introduction of statistical moments as a measure to quantify financial returns and risk has been recorded as one of the biggest milestones in modern portfolio theory. Markowitz (1952) was the first to offer a systematic approach to the problem of decision-making in the field of portfolio selection. He introduced a mean-variance criteria for identifying the most effective portfolio. In this framework, an optimal portfolio is achieved by maximizing expected returns represented by the first-order moment (mean) of the return distribution of the stocks, while minimizing the risk quantified by the second-order moment about the mean (variance). However, since variance works only for symmetric distributions, relying on it as a risk measure can sometimes lead to misleading results. To address this issue, Philippatos and Wilson (1972) introduced the mean-entropy framework, establishing entropy as a more general and competent risk measure due to its independence from distributional symmetry and non-reliance on numerical data. The mean-entropy model is shown to be equally efficient as the mean-variance model. Thereafter, entropy-based approaches have gained prominence in portfolio selection as alternatives to traditional risk measures. Philippatos and Gressis (1975) further validated this across uniform, normal and lognormal return distributions, emphasizing the advantages of entropy in not requiring distributional assumptions. Furthermore, Usta and Kantar (2011) and Aksarayli and Pala (2018) incorporated higher order statistical moments in their entropy-based models to develop the mean-variance-skewness-entropy (MVSE) and the mean-variance-skewness-kurtosis-entropy (MVSKE) models, respectively. Recently, Mercurio et al. (2020) introduced a return-entropy-based model that outperformed existing portfolio optimization approaches. Collectively, these studies indicate that integrating entropy with other decision-making criteria improves model performance by enhancing robustness in decision-making, especially in highly uncertain environments. (cf. Yang et al., 2017; Brito, 2023; and references therein).  \\
\hspace*{0.2in} However, all these studies just focused on constructing efficient portfolios from a pre-selected set of stocks. But real-world investment begins with identifying the ideal set of stock options. Increasing the number of stocks in a portfolio helps in lowering the variance of a portfolio. This will lead to reduction in risk. Consequently, diversifying a portfolio lowers the perceived risk associated with return uncertainty. However, the degree of diversification of a portfolio depends on the investors' risk attitudes towards return uncertainty and the expected utility of the stocks. It is considered that investors with high risk tolerance will prefer actions or investment decisions with higher returns while those with low risk tolerance will prefer decisions with moderate returns so as to avoid risk. So both risk aversion tendency based on the utility of the outcome and sensitivity towards return uncertainty needs to be captured.
Thus, the balance between risk and return can be achieved by monitoring the number of stocks prior to constructing an efficient frontier for portfolios. 
A variety of stock selection models are developed in literature. Among them some noteworthy contributions include the use of some machine learning methods designed to help with stock selection and construction of efficient portfolios. Liu and Yeh (2017) used artificial neural networks to develop stock selection models tailored to investors' risk preferences. Various other machine learning methods, such as long short-term memory (LSTM), eXtreme Gradient Boosting (XGBoost), support vector regression, genetic algorithms, random forests, and support vector machines (SVM), have also been applied to evaluate stock quality and rank them based on predicted risks and returns (cf. Brito, 2023, and references therein). These ranked stocks are then used to construct optimal portfolios using techniques discussed earlier, like the mean-variance and mean-entropy methods.
Furthermore, Yang et al. (2017) utilized the expected utility–entropy model, originally developed by Yang and Qiu (2005), for stock selection prior to portfolio optimization. It incorporated both the decision-maker's subjective preferences through a utility function and objective uncertainty quantified by Shannon entropy. They used the weighted linear average of the expected utility function and the entropy function to quantify the risk associated with an action. Yet, the EU-E model lacked mathematical generality and failed to satisfy certain behavioral axioms related to the ordering of preferences among uncertain prospects and their combinations. More recently, Brito (2020, 2023) extended this model by incorporating variance into the decision-making criteria for stock selection. These models have also helped in resolving several famous decision paradoxes and investment problems, offering a better explanation for anomalies in individual decisions. They also provided insights to the mean-variance criterion introduced by Pollatsek and Tversky (1970) and Markowitz (1952).\\
\hspace*{0.2in} But these decision analysis models are mostly limited to Shannon entropy, without exploring any fractional order generalizations of Shannon entropy or real-time adaptability. Nevertheless, the existing decision models fall short in capturing the dynamic nature of investors' risk tolerance levels in relation to the uncertainty embedded in stock return distributions. They also lack the integration of machine learning techniques, such as artificial neural networks (ANNs), for analyzing the non-linear relationships between the input and response variables. Besides, the choice of utility function and trade-off parameter also lacks empirical calibration or behavioral justification. Furthermore, the existing entropy-based stock selection models by Yang et al. (2017) and Brito (2023) fail to effectively capture the relationship between investors' risk sensitivity to return uncertainty and the expected utility of stocks.
To address the gaps in the existing literature, this paper presents a novel fractional order based entropy approach of representing investor preferences and risk tolerance in decision-making models under risk. The sensitivity of individuals towards risk due to return uncertainty is represented by the fractional entropy through means of the fractional parameter.\\
\hspace*{0.2in} Therefore, the rest of the paper is organized into the following sections: In Section 2, we have rigorously analyzed the behavior of the fractional order information gain and entropy functions due to Ubriaco with respect to the differences in individual risk attitudes in the human decision-making context. The results are further illustrated by real-time portfolio data.  
In Section 3, we introduce two new risk measures based on the Ubriaco fractional entropy for modeling decision-making under risk, examine their properties as risk measures and analyze their performance in resolving some important investment problems and decision paradoxes for different values of $ \alpha $ in comparison to existing Shannon entropy-based decision models. Section 4 focuses on introducing stock selection methods based on the proposed decision models using the PSI20 data and constructing efficient portfolio frontiers using mean-variance optimization technique. Finally, in Section 5, we conclude the present study. Throughout the paper, we have taken $ \log(\cdot) $ as the natural logarithm function.

\section{Information gain and entropy functions}
The Shannon entropy is one of the most popular information measure quantifying the uncertainty of a random action having multiple possible states of nature or outcomes pioneered by Shannon (1948). The uncertainty associated with a specific state can be represented in terms of information gain function. Then the entropy can be considered as the expected value of the information gain. High values of information gain (or the level of surprise) on the actual occurrence of an outcome is associated with higher uncertainty. The Shannon information gain function is defined as:
\begin{equation}
	I_S(p) = -\ln p(x),
\end{equation} where $ p(x) = P(X=x) $ is the probability mass function (pmf) describing the distribution of a discrete random variable $ X $ corresponding to an action. \\
\hspace*{0.2in} In the context of human decision-making, individuals’ differing levels of risk tolerance have a significant impact on the perceived uncertainty of an outcome. Risk-seeking individuals are more affirmative about the occurrence of low probable outcomes. Such people have an adventurous attitude while making decisions. And the ones who are more averse to risks are more doubtful about low probable outcomes. Such individuals are said to have a conservative attitude towards certainty of an outcome (see Aggarwal, 2021a;b). Unfortunately, this subjective uncertainty of outcomes depending on risk attitudes of individuals cannot be portrayed by the Shannon entropy. 
To remove this potential impediment of the Shannon entropy, Aggarwal (2021b) formulated a new form of probabilistic entropy that offers a holistic representation of the perceived uncertainty.
This formulation introduces a positive real-valued parameter to capture an individual's attitude towards uncertainty. Despite all the modifications in the entropy formulation, the interval of entropy values between minimum and maximum information gain is too wide for most of the existing generalized entropy formulations. Intuitively, this will affect the prediction accuracy and consistency of the decision models. Hence, the Shannon entropy is predominantly used in the decision models under risk and uncertainty. \\
\hspace*{0.2in} 
However, this issue can be handled effectively by the fractional order entropy (\ref{1.1}) with the help of the fractional parameter $ \alpha $ since $ \alpha $ is known to lie between 0 and 1. This influenced us to adopt the fractional entropy due to Ubriaco in modeling risky decisions.
The fractional entropy proposed by Ubriaco (2009) is particularly noteworthy, as it serves as a generalized fractional form of Shannon entropy. This generalization, featuring a fractional parameter $ \alpha \in (0,1) $, offers greater flexibility in describing complex systems with varying statistical characteristics with more granularity (see Machado, 2014;2020). Furthermore, its fractional formulation has been utilized to derive a significant class of fractional-order entropies by incorporating the cumulative distribution function or survival function instead of the probability mass function (see Di Crescenzo et al., 2021 and Xiong et al., 2019). If $ p(x) = P(X=x) $ is the probability mass function for $ x^{th} $ outcome of an action or event, then the fractional order entropy for a $n$-state discrete system is defined by Ubriaco as:
\begin{equation}
	H_U^\alpha (p) = \sum_{x=1}^{n}p(x) (-\log p(x))^\alpha;~ \alpha \in [0,1].
	\label{1.1}
\end{equation} 
Clearly, when $ \alpha = 1 $, the fractional order entropy (\ref{1.1}) is equivalent to the Shannon entropy given by 
\begin{equation}
	H_S(p) = -\sum_{x=1}^{n}p(x) \log p(x).
	\label{1.2}
\end{equation}
Here, the parameter $ \alpha $ allows the Ubriaco entropy function (\ref{1.1}) to customize itself so as to accommodate the unique attitudes of decision-makers towards the likelihood of an outcome.
\subsection{Ubriaco information gain function}
In decision-making, individuals exhibit either adventurous or conservative attitudes, reflecting high or low risk tolerance, respectively. An adventurous person tends to trust the occurrence of an outcome even when it is uncertain or unlikely. In contrast, a conservative individual remains doubtful about an outcome’s occurrence, even as its probability increases. To capture these nuances in human behavior, we use the Ubriaco information gain function, expressed as follows:
\begin{equation}
I_U^\alpha(p) = (-\log p(x))^\alpha, \alpha \in [0,1],
\label{2.3}
\end{equation}where $ p(x) \in (0,1]$ denotes probability of the $ x^{th} $ likely outcome of a random event.\\
\hspace*{0.2in} From (\ref{2.3}), we observe that the information gained from an outcome is inversely proportional to its likelihood $ p(x). $ However, the rate at which information decreases with increasing $ p(x) $ varies with different values of $ \alpha. $ Therefore, it is important to examine the behavior of function (\ref{2.3}) with respect to both varying probabilities and changing $ \alpha $ values. \\
\begin{figure}[t]
\includegraphics[width=0.45\textwidth]{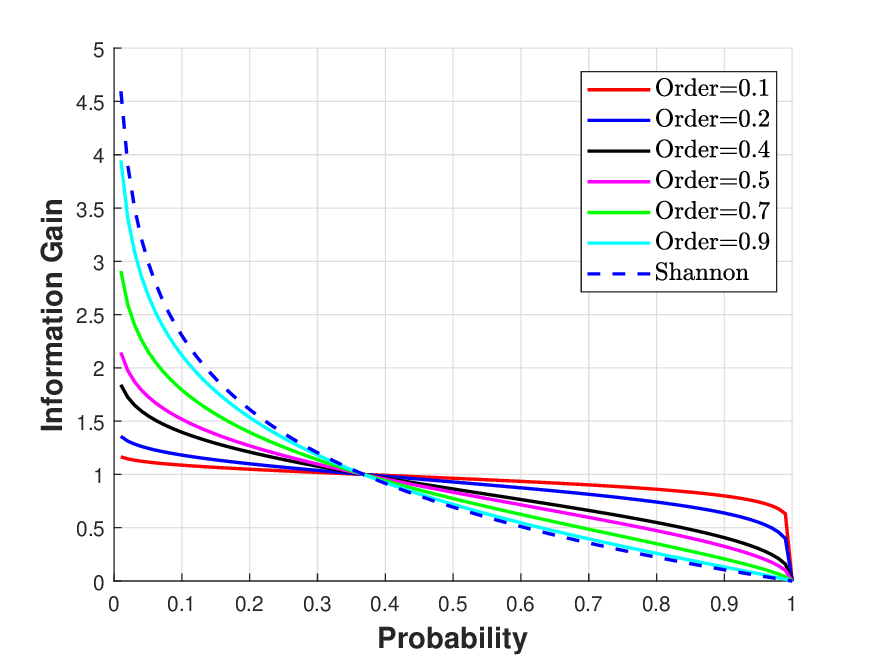}~~~ \includegraphics[width=0.45\textwidth]{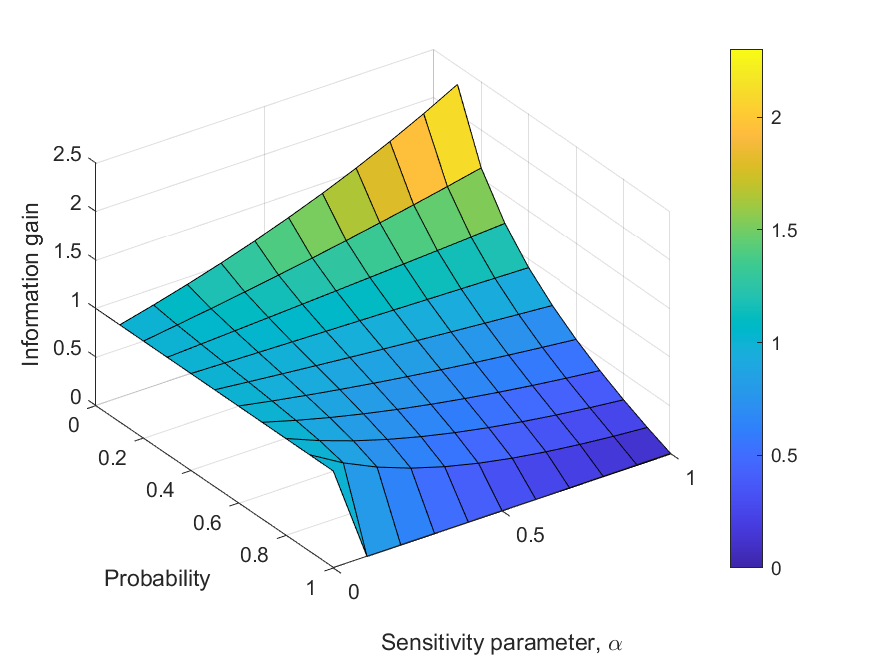}
\caption{Information gain function for Ubriaco entropy w.r.t probability.}
\label{info_gain}
\end{figure}
\hspace*{0.2in} Fig. \ref{info_gain} shows that for higher values of $ \alpha \geq 0.5 $, the sharper decline in information gain at low $ p $ reflects a more adventurous attitude. Such individuals gain significant information from unlikely events but become nearly certain as $ p $ approaches 1, leading to minimal change in the perceived uncertainty of an outcome.
Conversely, when $ \alpha < 0.5 $, the curve flattens, indicating a slower decline in information gain, especially at lower probabilities ($ p < 0.5 $). As $ p $ increases beyond 0.5, the information gain decreases more rapidly, dropping sharply to zero near $ p = 1 $, where the decision-maker becomes fully confident. At this point, the individual experiences little to no surprise or information gain on  the actual occurrence of the outcome. Thus, lower $ \alpha $ values reflect a more conservative attitude toward uncertainty.
Therefore, as an individual's adventurous nature increases, the subjective range of low-probability values narrows with rising $ \alpha $ (from 0.5 to 1), and widens as the conservative attitude dominates. This suggests that higher $ \alpha $ values exert a stronger influence on reducing the gain function at low probabilities, while lower $ \alpha $ values (near 0) have a greater effect at higher probabilities $ p $. \\
\hspace*{0.2in} Thus, the deviations in individual risk-taking tendencies, unaccounted for by the Shannon entropy, can be effectively captured using the Ubriaco information gain function (\ref{2.3}). By adjusting $ \alpha $ between 0 and 1, the function reflects how individuals with different risk tolerance levels perceive information gain across varying probability levels. Values of $ \alpha $ near 1 correspond to highly adventurous attitudes, those near 0 to highly conservative ones, and values around 0.5 indicate an almost risk-neutral stance. This conservative behavior of an individual under uncertainty, overlooked by Shannon entropy, is well-represented through fractional-order entropy.

\subsection{Sensitivity of Ubriaco information gain function}
Since variations in human nature significantly influence the rate at which information about the certainty of an outcome decreases, the responsiveness of the information gain function to changes in probability exhibits distinct patterns for different values of $ \alpha $, particularly at low and high values of $ p $. This responsiveness is referred to as the elasticity or sensitivity of the information gain function. \\
\hspace*{0.2in}The sensitivity can be quantified by taking the modulus of the ratio of the function's slope with respect to probability to the information gain per unit probability. Mathematically, it is expressed as:
\begin{equation}
\mathcal{E} = \left|\frac{dI_U^\alpha(p)/dp}{I_U^\alpha(p)/p}\right| = \alpha/(-\log p(i)),
\label{2.4}
\end{equation} where $|\cdot|$ represents the absolute value. 
$ \mathcal{E} $ indicates the degree of the information gain function's responsiveness to variations in likelihood values. \\ 
\begin{figure}[t]
\centering
\includegraphics[width=0.45\textwidth]{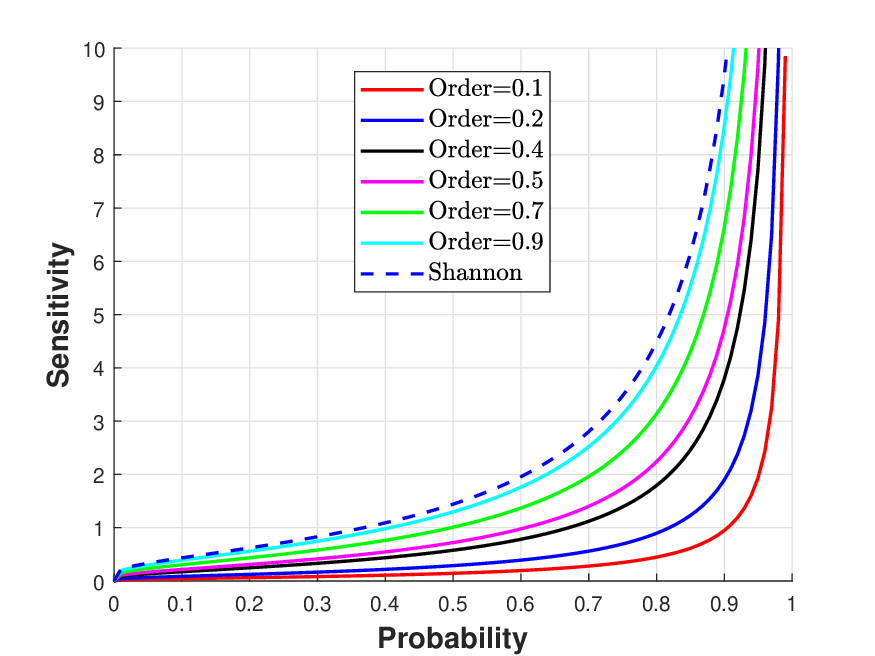} ~~~ \includegraphics[width=0.45\textwidth]{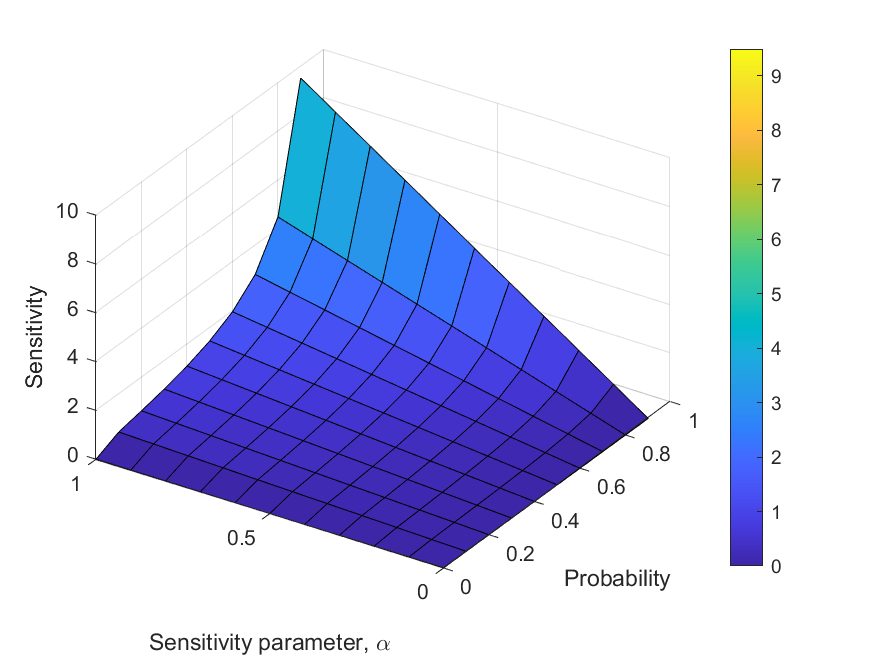}
\caption{Sensitivity of information gain function for Ubriaco entropy.}
\label{elasticity}
\end{figure}
\hspace*{0.2in} From the sensitivity plot in Fig. \ref{elasticity}, the sensitivity function (\ref{2.4}) appears as an increasing and convex function for all $ \alpha \in [0,1] $. Information gain is more sensitive to higher probability values ( $ p > 0.5 $) and less responsive at lower $ p $. As $ \alpha $ decreases, the values of $ p $ where sensitivity sharply increases shift closer to 1. This indicates that lower $ \alpha $ values correspond to reduced responsiveness across a wider range of probabilities, reflecting a more conservative attitude. For smaller $ \alpha $, significant changes in information gain occur only when $ p $ approaches 1. In other words, as $ \alpha $ nears 0, noticeable changes in sensitivity occur within a narrower band of high probability values (with $ p $ shifting toward 1). \\
\hspace*{0.2in} Thus, the fractional order $ \alpha $ can be tuned to represent varying sensitivity levels aligned with individuals' attitudes towards uncertainty. Clearly, $ \alpha $ acts as a sensitivity parameter for the Ubriaco gain function, effectively capturing variations in behavior, especially among low-risk takers. 

\subsection{Ubriaco entropy function}
Entropy can be interpreted as the area under the curve of an entropy function. Comparing these areas helps evaluate the effectiveness of an entropy measure in representing uncertainty. Intuitively, an entropy function that covers a larger area is considered a more efficient measure of information for a given random variable with finite support (see Karci, 2016). Therefore, we analyze the Ubriaco entropy function \( h_U^\alpha(p) = p(x)[-\log(p(x))]^\alpha \) for different values of \( \alpha \) and compare its performance with the Shannon entropy \( h_S(p) = -p(x)\log p(x) \). The following theorem presents the monotonic behavior of \( h_U^\alpha(p) \) as the order \( \alpha \) varies. The proof is straightforward and hence omitted.

\begin{t1}
The entropy function $ h_U^\alpha $ decreases as $ \alpha $ decreases from 1 to 0 when $ 0 < p(x) \leq \frac{1}{e} (= 0.368) $.
However, when $ 0.368 < p(x) < 1 $, the entropy function shows an increasing behavior with decrease in $ \alpha $. 
\label{thm 2.1}
\end{t1}

\begin{r1}
From the above theorem, we can conclude that for the case when $ \frac{1}{e} < p(x) < 1 $, the fractional order entropy measures with order $ \alpha < 1 $ are preferred to represent the uncertainty of a system as it contains more information than the Shannon entropy.
Nevertheless, when $ 0 < p(x) \leq \frac{1}{e} $, the Shannon entropy will dominate all other cases of Ubriaco entropy with order $ \alpha < 1 $.  \\
\hspace*{0.2in}This behavior can be understood better through the illustration presented in Fig. \ref{ent_order} for a $n$-state system. It also shows that the entropy function $  h_U^\alpha(p) $ is non-negative and concave for all values of $ p. $ 		
\begin{figure}[h]
\centering
\includegraphics[width=0.5\textwidth]{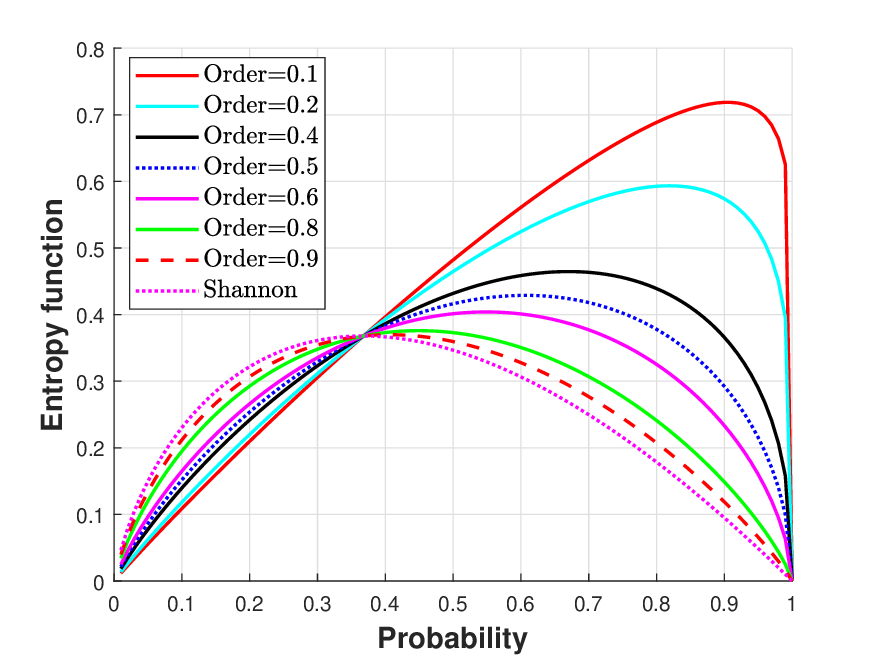} 
\caption{Ubriaco entropy function of a $ n $ -state system for different $ \alpha. $}
\label{ent_order}
\end{figure} 
\label{rem 2.1} 
\end{r1}
Thus, for events with outcome probabilities greater than 0.368, the parameter \( \alpha \) in the fractional-order entropy can be tuned to extract maximum information about uncertainty, aiding in more informed decision-making. 
\begin{figure}[h]
\centering
\includegraphics[width=0.5\textwidth]{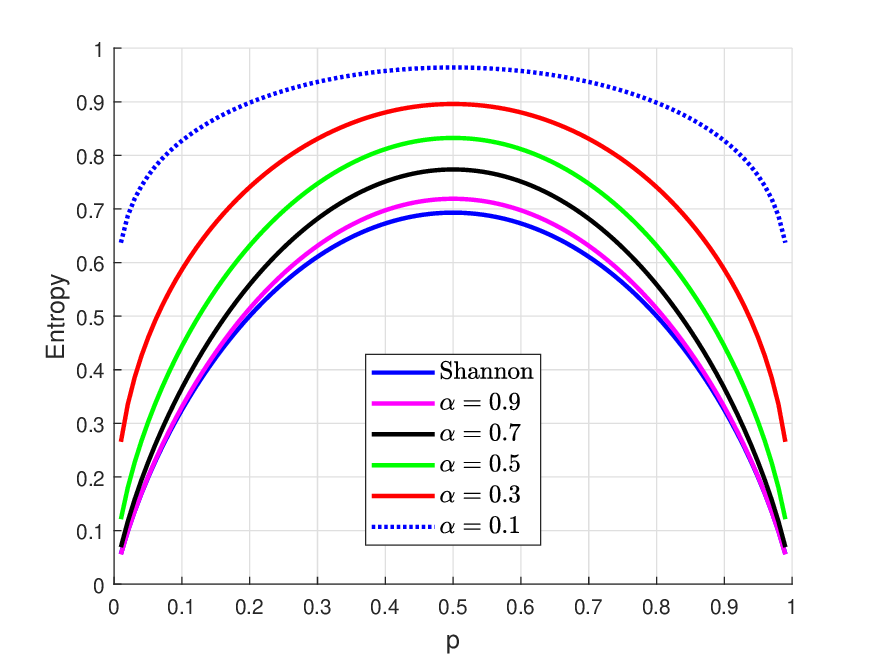} 
\caption{Ubriaco entropy of a binomial system for different $ \alpha. $}
\label{Fig6} 
\end{figure} 	
This claim is further supported by analyzing the behavior of Ubriaco entropy across varying probabilities and \( \alpha \) values for a binomial system with two microscopic states, as illustrated in Fig.~\ref{Fig6}. For such a system, the fractional-order entropy is given by:
\begin{equation}
H_U^\alpha(p) = p[-\log p]^\alpha + q[-\log q]^\alpha,
\label{3.17}
\end{equation}where $ p = $probability of success of an outcome and $ q = 1-p $ is the probability of its failure. \\
\hspace*{0.2in} From Fig.~\ref{Fig6}, it is evident that the fractional-order entropy yields higher values, i.e., more information about uncertainty, when \( \alpha = 0.1 \), and the least information when \( \alpha = 1 \) (corresponding to Shannon entropy). This indicates that the information content of Ubriaco entropy increases as \( \alpha \) decreases in a binomial system. The figure also shows that, for all \( \alpha \in (0,1] \), Ubriaco entropy attains its peak when the outcomes are equally likely, i.e., when \( p = 0.5 \).

\subsection{Validation through portfolio data}
\begin{figure}[h]
\centering
\includegraphics[width=0.45\textwidth]{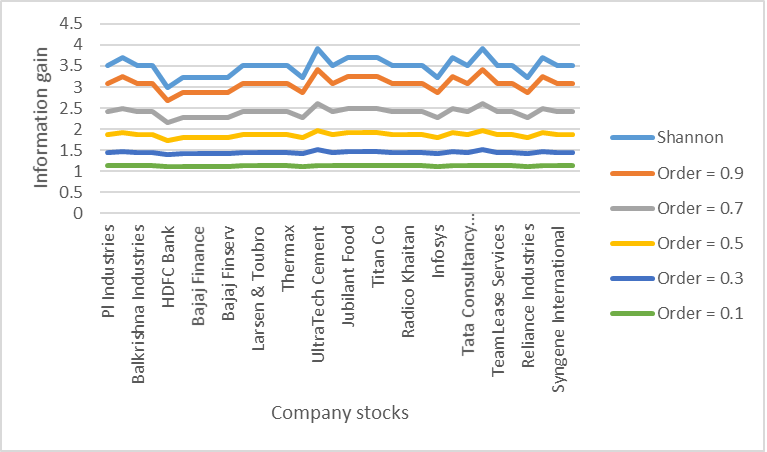} ~~~ \includegraphics[width=0.45\textwidth]{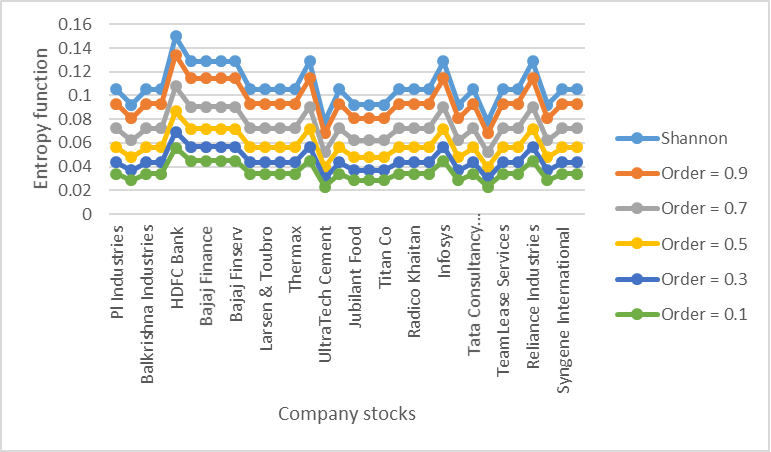}
\caption{Variation of information gain (left) and entropy (right) functions with different $ \alpha $ values for the Di portfolio.}
\label{fig:div}
\end{figure}

\begin{figure}[h]
\centering
\includegraphics[width=0.45\textwidth]{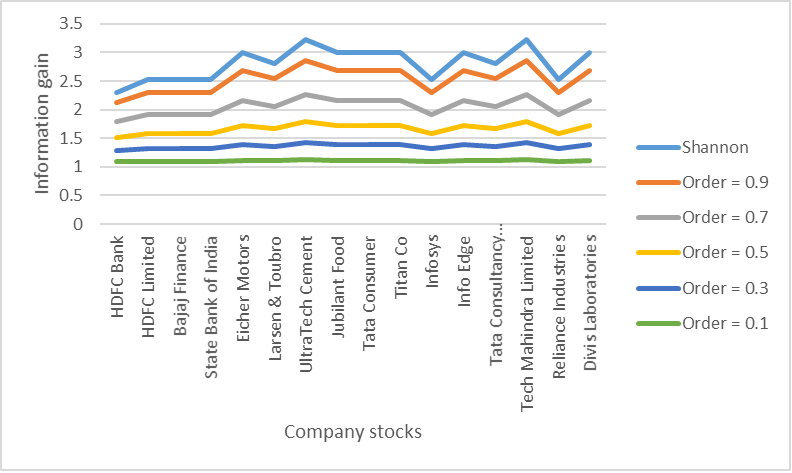} ~~~ \includegraphics[width=0.45\textwidth]{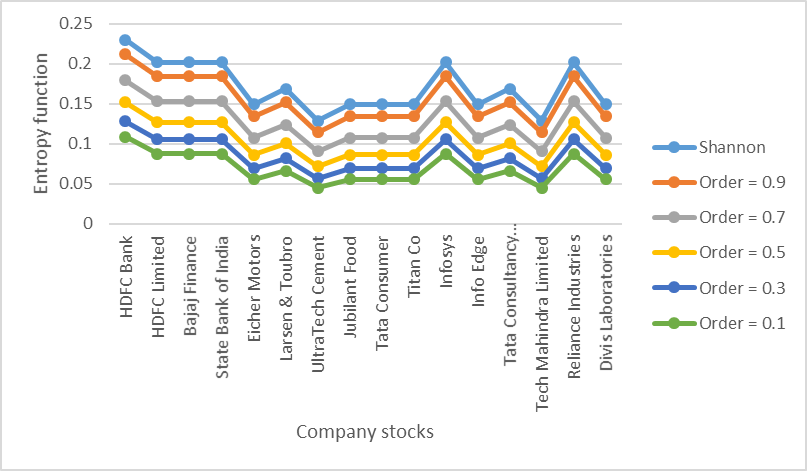}
\caption{Variation of information gain (left) and entropy (right) functions with different $ \alpha $ values for the Lc portfolio.}
\label{fig:lc}
\end{figure}

\begin{figure}[h]
\centering
\includegraphics[width=0.45\textwidth]{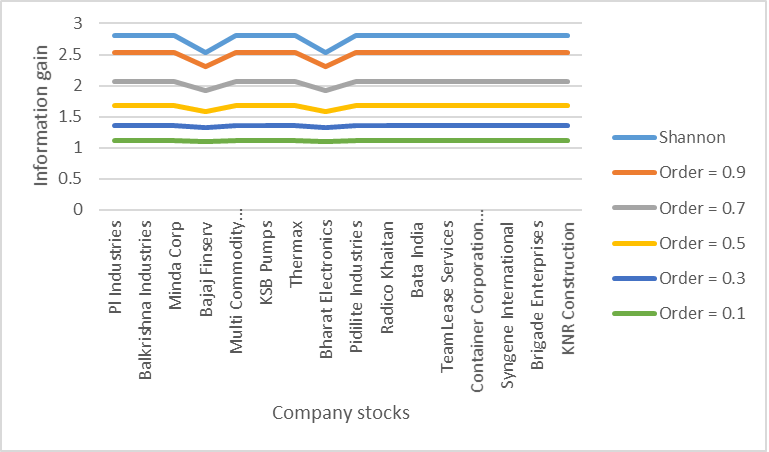} ~~~  \includegraphics[width=0.45\textwidth]{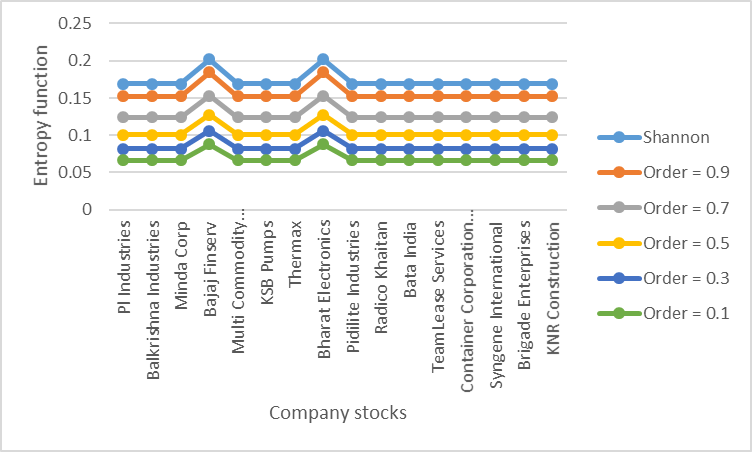}
\caption{Variation of information gain (left) and entropy (right) functions with different $ \alpha $ values for the Mc portfolio.}
\label{fig:mc}
\end{figure}

\begin{figure}[h]
\centering
\includegraphics[width=0.45\textwidth]{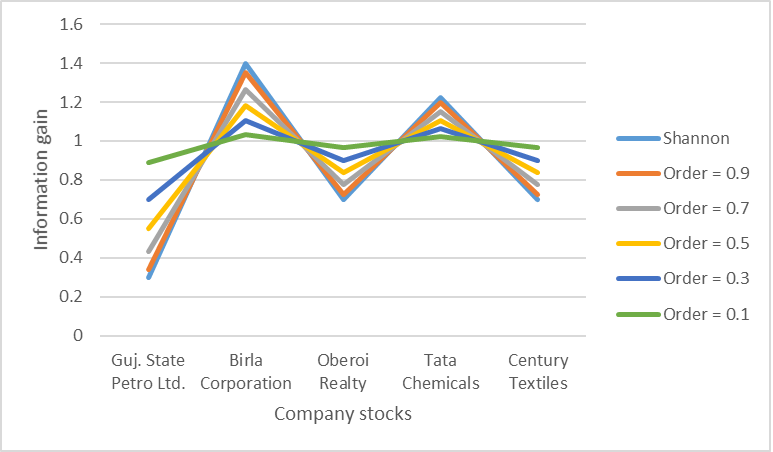} ~~~ \includegraphics[width=0.45\textwidth]{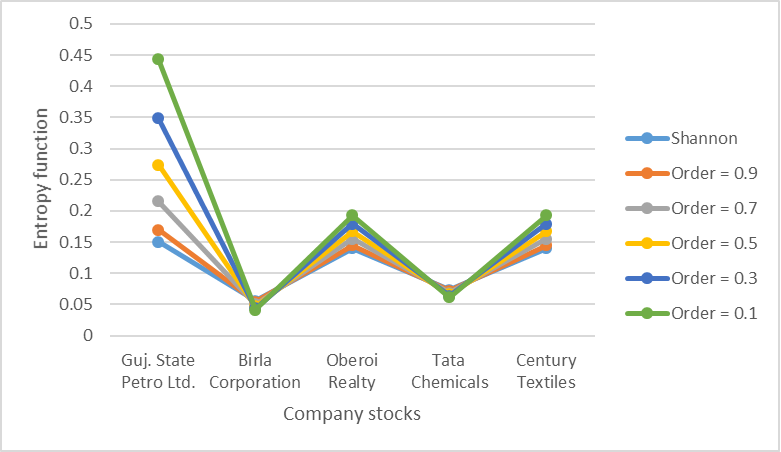}
\caption{Variation of information gain (left) and entropy (right) functions with different $ \alpha $ values for the Hy portfolio.}
\label{fig:Hyp}
\end{figure}

\noindent
The previous section provided a rigorous analysis of the monotonicity of Ubriaco information gain function with respect to probabilities through graphical representations. It was observed that for individuals with higher risk tolerance (\( \alpha > 0.5 \)), the effect of \( \alpha \) is more pronounced in reducing uncertainty for low-probability outcomes as their likelihood increases. Consequently, entropy reaches its maximum when \( \alpha \to 1 \) at low probability values (specifically, below 0.368, as stated in Theorem~\ref{thm 2.1}).
In this section, we examine four portfolio models, namely, diversified, large-cap, mid-cap, and hypothetical, shown in Tables~\ref{tab:di}, \ref{tab:lc}, \ref{tab:mc}, and \ref{tab:hy}, respectively, to validate our theoretical findings on how fractional-order Ubriaco entropy behaves with varying risk attitudes. The data sets used in this section, displayed in Tables \ref{tab:di}–\ref{tab:mc} are collected directly from an actual investor interested to invest in the Indian stock market. The entropy values for each portfolio are computed using the Ubriaco formula~(\ref{1.1}) and summarized in Table~\ref{tab:entropy}. The corresponding ordering of entropy values with increasing \( \alpha \) is presented in Table~\ref{tab:order1}. \\
\hspace*{0.2in} For the first three real-time portfolio models (Tables~\ref{tab:di}–\ref{tab:mc}), the entropy values increase with \( \alpha \), yielding the maximum value for Shannon entropy (\( \alpha = 1 \)), consistent with theoretical expectations. However, the hypothetical portfolio (Table~\ref{tab:hy}) exhibits a contrasting trend, where entropy is highest for \( \alpha = 0.1 \) and lowest for Shannon entropy. This reversal in the trend indicates that entropy functions with higher \( \alpha \) values give less information gain for high-return outcomes in this case, while lower \( \alpha \) values capture more uncertainty. Such behavior, distinguishable from the other portfolios, aligns with Theorem~\ref{thm 2.1}.
Thus, these four types of portfolio help in examining the sensitivity of entropy to changes in investors’ risk tolerance (\( \alpha \)) and varying combination of weights (represented by $ p $ values) for portfolio stocks.
Figs.~\ref{fig:div}–\ref{fig:Hyp} visually depict the variation of the information gain function \( I_U^\alpha(p) \) and entropy function \( h_U^\alpha(p) \) across portfolio categories as \( \alpha \) changes. In particular, when \( \alpha \) is close to zero, both functions remain nearly unchanged with respect to \( p(i) \), highlighting increased conservatism at lower \( \alpha \) values.  \\
\hspace*{0.2in} Furthermore, Figs.~\ref{fig:div}, \ref{fig:lc}, \ref{fig:mc}, and \ref{fig:Hyp} visually illustrate the behavior of the information gain function \( I_U^\alpha(p) \) and entropy function \( h_U^\alpha(p) \) across different portfolio categories as \( \alpha \) varies. Figs. \ref{fig:div}, \ref{fig:lc} and \ref{fig:mc} demonstrate that the information gain function remains nearly constant when both \( p \) and \( \alpha \) are low ($< 0.5$), indicating increased conservatism as \( \alpha \) approaches zero. However, relatively higher fluctuations are observed with the gain functions in the case of the hypothetical portfolio set as shown in Fig. \ref{fig:Hyp} for all the selected values of $ \alpha $, indicating higher sensitivity of the gain function at comparatively higher $ p $ values ($ p\geq 0.5 $) represented by the stock weights $p(x)$ of the hypothetical portfolio (cf. Table \ref{tab:hy}). Further the behavior of the corresponding  entropy functions illustrated by Figs. \ref{fig:div} - \ref{fig:Hyp} strongly support our findings stated in Theorem \ref{thm 2.1} and Remark \ref{rem 2.1}.\\
\hspace*{0.2in} Thus, these portfolio datasets provide an empirical evidence of the sensitivity of the information gain and entropy functions with changes in $ \alpha $ and $ p $ values. This reinforce the relevance of fractional-order entropy in modeling human decision-making events such as portfolio stock selection. By capturing subjective uncertainty based on individual risk tolerance, the Ubriaco entropy proves to be a valuable tool for decision analysis under risk and uncertainty. 

\section{Fractional entropy-based decision analysis models}
In the context of decision-making among risky choices, two key factors need to be considered: (1) greater risk is associated with increased uncertainty in the outcomes due to the occurrence of different states, and (2) lower risk is linked to outcomes that have a higher expected utility as perceived by decision-makers (see Yang and Qiu, 2005 and the references therein).\\
Bearing these factors in mind, Yang and Qiu (2005) expanded the classical decision analysis model under risk to better capture real-world scenarios where the distribution of states of nature varies across different elements of the action space. They developed a new decision model, called expected utility-entropy (EU-E) model, that integrates the expected utility of an action with the Shannon entropy of the state of nature through a risk tradeoff factor. This EU-E model demonstrated consistency with previous empirical findings on individual risk preferences. Further, Brito (2020) enhanced the EU-E model by incorporating variance known as expected utility-entropy variance (EU-EV) model as an additional risk measure to address certain decision paradoxes where the original EU-E model fell short. \\
\hspace*{0.2in} The general decision analysis model is framed as:
Let \( K \) be a finite action space comprising actions \( k_i \) for \( i = 1, 2, \dots, l \). Each risky action \( k_i \) is associated with a state space 
\[
\Theta_i = \{\theta_{i_1}, \theta_{i_2}, \dots, \theta_{i_{l_i}}\}.
\]
Let the overall state space be denoted by 
\[
\Theta = \{\theta_1, \theta_2, \dots, \theta_l\}.
\]
The probability distribution over states given action \( k_i \) is represented by 
\[
p_{ij} = P(\theta = \theta_{ij} \mid k = k_i),
\]
subject to the conditions \( \sum_{j=1}^{l_i} p_{ij} = 1 \) and \( p_{ij} \geq 0 \).
The expected payoff resulting from state \( \theta_{ij} \) under action \( a_i \) is given by 
\[
X = X(k_i, \theta_{ij}) = x_{ij}, \quad (i = 1, 2, \dots, l; \; j = 1, 2, \dots, l_i).
\]
The utility function, capturing an individual's risk attitude, is denoted by \( u = u(X) \). A detailed formulation of this generalized decision model can be found in Table 1 of Yang and Qiu (2005).
To advance the decision models proposed by Yang and Qiu (2005) and Brito (2020), we introduce the use of Ubriaco fractional entropy~(\ref{1.1}) in place of Shannon entropy to capture the perceived uncertainty of outcomes due to varying risk attitudes of individuals alongside the risk preferences or risk aversion tendency represented by the expected utility of the outcome. This leads us to formulate two novel risk measures: (1) the expected utility-fractional entropy (EU-FE) measure of risk, and (2) the expected utility-fractional entropy-variance (EU-FEV) measure of risk.
\begin{d1}
Consider a general decision analysis model $ G = (\Theta,A,u) $, with an increasing utility function $ u = u(X(a,\theta)) $ corresponding to action $ a \in A. $ Then, the expected utility- fractional entropy (EU-FE) measure of risk associated with an action $ a $ is defined by
\begin{equation}
R(a) = \lambda H_a^\alpha(\theta)-(1-\lambda)\frac{E[u(X(a,\theta))]}{max_{a\in A}\{|E[u(X(a,\theta))]|\}},
\label{3.8}
\end{equation}
where $ \lambda \in \mathcal{R} $ is the risk-tradeoff factor such that $ 0 \leq \lambda \leq 1 $ and $ H_a^\alpha(\theta) $ is the fractional entropy of the distribution of the state of nature corresponding to action $ a. $ Here, we assume that $ max_{a\in A}\{|E[u(X(a,\theta))]|\} $ exists and is non-zero. When $ max_{a\in A}\{|E[u(X(a,\theta))]|\} = 0, $ then $$ R(a) = \lambda H_a^\alpha(\theta). $$
\label{def_3.1}	
\end{d1}
\vspace{-3.0em}
\begin{d1}
Given a general decision analysis model $ G = (\Theta,\mathcal{B},u) $, with an increasing utility function $ u = u(X(b,\theta)) $ corresponding to action $ b \in \mathcal{B}. $ Then, the expected utility- fractional entropy and variance (EU-FEV) measure of risk associated with an action $ b $ is defined by
\begin{equation}
R(b) = \frac{\lambda}{2}\Bigg[H_b^\alpha(\theta) + \frac{Var[X(b,\theta)]}{max_{b\in \mathcal{B}}\{Var[X(b,\theta)]\}}\Bigg]-(1-\lambda)\frac{E[u(X(b,\theta))]}{max_{b\in \mathcal{B}}\{|E[u(X(b,\theta))]|\}},
\label{3.9}
\end{equation}
where $ \lambda \in \mathcal{R} $ is the risk-tradeoff factor such that $ 0 \leq \lambda \leq 1 $ and $ H_b^\alpha(\theta) $ is the fractional entropy of the distribution of the state of nature corresponding to action $ b. $ Here, we assume that $ max_{b\in \mathcal{B}}\{|E[u(X(b,\theta))]|\} $ exists and is non-zero. When $ max_{b\in \mathcal{B}}\{|E[u(X(b,\theta))]|\} = 0, $ then $$ R(b) = \frac{\lambda}{2}\Bigg[H_b^\alpha(\theta) + \frac{Var[X(b,\theta)]}{max_{b\in \mathcal{B}}\{Var[X(b,\theta)]\}}\Bigg]. $$
\label{def_3.2}
\end{d1}
\vspace{-1.5em}
\hspace*{0.2in} To assess decisions made by individuals across a finite set of actions based on a risk measure, the risk values ($ R(y), y= a ~\text{or}~ b $) associated with each action are compared. These risk values are determined by the subjective uncertainty of the outcomes and their perceived utility. The action with the lowest $ R(y) $  value is selected. Therefore, the EU-FE and EU-FEV decision models can be defined as:
\begin{d1}
Let $ G = (\Theta,Y,u) $ be a decision analysis model.
\begin{enumerate}
\item Let the action space $ Y = A ~\text{or}~ \mathcal{B} $ consist of two actions $ y_1, y_2 \in Y $ with related risk measures $ R(y_1) $ and $ R(y_2), $ respectively. Then the following results follow: \\ 
$ (i) $ action $ y_1 $ is selected, i.e., $ y_1 > y_2 $, if $ R(y_1) < R(y_2); $ \\
$ (ii) $ action $ y_2 $ is selected, i.e., $ y_1 < y_2 $, if $ R(y_1) > R(y_2); $ \\
$ (iii) $ no action is superior i.e., any action can be selected or, $ y_1 = y_2 $, if $ R(y_1) = R(y_2). $ 
\item If $ Y $ is a discrete set of $ n $ objects, i.e., $ Y = \{y_1,y_2,\dots,y_n\}, $ then the EU-FE and the EU-FEV order of preferences of actions will be in accordance with decreasing order of risk values. The most preferred action will then be 
\[R(y_i) = min_{y_j\in Y} R(y_j). \] 
\label{def_3.3}
\end{enumerate}
\end{d1} 
\vspace{-1.5em}
\hspace*{0.2in} We can observe that when $ \alpha = 1, $ both the EU-FE and EU-FEV decision models reduce to the EU-E and EU-EV models, respectively. Moreover, for $ \lambda = 0, $ the proposed decision models reduce to the expected utility model ($ y_1 > y_2 $ or $ R(y_1) < R(y_2) $ iff $E[u(y_1)] > E[u(y_2)] $). 
\subsection{Properties of the EU-FE and EU-FEV measures of risk}
The EU-FE and EU-FEV risk measures adhere to the core principles of risk measures, similar to those of previously established risk measures. This section will explore some of the key properties of these proposed risk measures.
\begin{t1}
Let $ a_1 $ and $ a_2 $ be two different actions in the action space associated with a general decision model $ G = (\Theta,A,u) $ having equal expected utility, $ H_{a_1}^\alpha(\theta) $ and $ H_{a_2}^\alpha(\theta) $ denote fractional entropies of the corresponding states of nature, then $ R(a_1) < R(a_2) $ when $ H_{a_1}^\alpha(\theta) < H_{a_2}^\alpha(\theta) $ for some $ \alpha. $ For the EU-FEV model, $ R(b_1) < R(b_2) $ if and only if $$ \Big[H_{b_1}^\alpha(\theta) + \frac{Var[X(b_1,\theta)]}{max_{b_1\in \mathcal{B}}\{Var[X(b_1,\theta)]\}}\Big] < \Big[H_{b_2}^\alpha(\theta) + \frac{Var[X(b_2,\theta)]}{max_{b_2\in \mathcal{B}}\{Var[X(b_2,\theta)]\}}\Big]. $$ 
\end{t1}
\textbf{\textit{Proof.}}
The result follows from Definitions \ref{def_3.1} and \ref{def_3.2}. $\hfill\square$\\

\hspace*{0.2in} Intuitively, according to the theorem mentioned above, an action with lower entropy and smaller variance will be associated with reduced risk.                       
\begin{t1}
Let $ G = (\Theta,A,u) $ be a general decision analysis model, with the action space comprised of two actions $ a_1 $ and $ a_2 $, $ H_{a_1}^\alpha(\theta) $ and $ H_{a_2}^\alpha(\theta) $ denoting the fractional entropies of the corresponding states of nature $ \theta $, respectively, then $ R(a_1) < R(a_2) $ when $ H_{a_1}^\alpha(\theta) = H_{a_2}^\alpha(\theta) $ for a particular value of $ \alpha $ and $ E[u(X(a_1,\theta))] > E[u(X(a_2,\theta))]. $ Further, when $ H_{b_1}^\alpha(\theta) = H_{b_2}^\alpha(\theta) = 0 $, then $ R(b_1) < R(b_2) $ if and only if $E[u(X(b_1,\theta))] > E[u(X(b_2,\theta))]$. 	
\label{thm_4.1}
\end{t1}
\textbf{\textit{Proof.}}
This result is also an immediate consequence of Definitions \ref{def_3.1} and \ref{def_3.2} and of the fact that for a decision model $ G(\Theta,Y,u), $ we have $ Var(X(y,\theta)) = 0 $ if $  H_{y}^\alpha(\theta) = 0 $. $\hfill\square$\\

\hspace*{0.2in} From Theorem \ref{thm_4.1}, we get that the action with higher expected utility will be less riskier when the fractional entropy values of the corresponding outcomes are equal. Further, the validity of the EU-FE and EU-FEV risk measures defined in (\ref{3.8}) and (\ref{3.9}) is examined by exploring several properties in relation to prior studies on decision models involving Shannon entropy and individual's perception about risk.  

\begin{t1}
Let $ G = (\Theta,Y,u) $ be a general decision analysis model with non-negative outcomes:
\begin{enumerate}
\item Let the action space be described as $ Y = \{y, y+k\}, $ where $ k > 0 $ is a constant, then $$ R(y + k) < R(y), $$  
i.e., risk decreases with the addition of a positive constant to all outcomes of an action.
\item If the action space contains two elements, namely, $ Y = \{y,ty\}, $ where $ t > 1 $ is a constant, then for the EU-FE decision model, $$ R(ty) < R(y), $$
i.e., EU-FE risk decreases if all the outcomes of an action are multiplied by a positive constant.\\
\hspace*{0.2in} For the EU-FEV decision model, $$ R(ty) < R(y), $$ if and only if \[\lambda \in \Bigg[0,~\frac{1-\frac{E[u(X)]}{E[u(kX)]}}{\frac{3}{2}-\frac{1}{2k^2}-\frac{E[u(X)]}{E[u(kX)]}}\Bigg).\]	\end{enumerate} 
\end{t1}
\textit{Proof:} \begin{enumerate}
	\item As we know that the utility function $ u(X) $ is increasing, then for $ k>0, $ we have $$ u(X+k) > u(X). $$ If $ Y = \{y,y+k\}, $ then $$ max_{y\in Y}\{|E[u(X(y,\theta))]|\} = E[u(X+k)] ~\text{and}~ Var(X+k)= Var(X). $$ Hence, from Definition \ref{def_3.2}, we get that 
	\[R(b) = \frac{\lambda}{2}\Big[H_b^\alpha(\theta) + 1\Big]-(1-\lambda)\frac{E[u(b)]}{E[u(b+k)]}\]
	and \[R(b+k) = \frac{\lambda}{2}\Big[H_{b+k}^\alpha(\theta) + 1\Big]-(1-\lambda).\]
	Similarly, from Definition \ref{def_3.1}, we have that
	\[R(a) = \lambda H_a^\alpha(\theta) - (1-\lambda)\frac{E[u(a)]}{E[u(a+k)]}\]
	and \[R(a+k) = \lambda H_{a+k}^\alpha(\theta)-(1-\lambda).\]
	Since, $ H_{y+k}^\alpha(\theta) = H_y^\alpha(\theta)$ ($ y = a ~\text{or}~ b $) for a fixed value of $ \alpha \in [0,1], $ one gets that $ R(y+k) < R(y). $ 	
	\item As we have an increasing utility function and non-negative outcomes $ X = X(b,\theta) $ such that $ E[X] = 0, $ then we get $ u(tX) > u(X) $ for $ t>1. $ Since, $ \mathcal{B} = \{b,tb\}, $ therefore
	\[max_{b\in \mathcal{B}}\{|E[u(X(b,\theta))]|\} = E[u(kX)] ~\text{and}~ Var(kX)= k^2E[X^2] > Var(X) = E[X^2].\]
	Hence, from Definition \ref{def_3.2}, we get that 
	\[R(b) = \frac{\lambda}{2}\Bigg[H_b^\alpha(\theta) + \frac{1}{k^2}\Bigg]-(1-\lambda)\frac{E[u(X)]}{E[u(kX)]}\]
	and \[R(kb) = \frac{\lambda}{2}\Big[H_{kb}^\alpha(\theta) + 1\Big]-(1-\lambda).\]
	Since entropy depends on the probabilities of the states of nature $ p(i) $ and $ \alpha, $ for a particular value of $ \alpha,~ H_{kb}^\alpha(\theta) = H_b^\alpha(\theta) $. Hence, it follows that $ R(kb) < R(b) $ if and only if
	\[\lambda\Bigg(\frac{3}{2}-\frac{1}{2k^2}-\frac{E[u(X)]}{E[u(kX)]}\Bigg)-\Bigg(1-\frac{E[u(X)]}{E[u(kX)]}\Bigg) > 0.\]
	Therefore, we have that for $ R(kb) < R(b) $ if and only if 
	\[\lambda \in \Bigg[0,~\frac{1-\frac{E[u(X)]}{E[u(kX)]}}{\frac{3}{2}-\frac{1}{2k^2}-\frac{E[u(X)]}{E[u(kX)]}}\Bigg).\]
	Similarly, it follows from Definition \ref{def_3.1} that
	\[R(a) = \lambda H_a^\alpha(\theta)-(1-\lambda)\frac{E[u(X)]}{E[u(kX)]}\]
	and \[R(ka) = \lambda H_{ka}^\alpha(\theta) -(1-\lambda).\]
	Thus, we obtain $ R(ka) < R(a) $ since for actions $ a $ and $ ka, ~ H_{ka}^\alpha(\theta) = H_a^\alpha(\theta). $ $\hfill\square$ 
\end{enumerate}

\subsection{Decision problems and solutions}
In order to investigate the effectiveness of the proposed risk measures, we consider the following decision problems and paradoxes used by Brito (2020). Additionally, by comparing these models to conventional decision-making approaches (EU-E and EU-EV), we highlight their distinct benefits in evaluating complex investment decisions and addressing well-known decision-making paradoxes.

\subsubsection{Nawrocki and Harding investment problem}
Nawrocki and Harding (1986) introduced a security investment problem, where two securities have the same Shannon entropy value but different risk levels. This showed the inadequacy of the Shannon entropy as a measure of security risk. So he defined a state-value weighted entropy to depict the dispersion of frequency classes of securities. The details of this problem can be found in Table \ref{tab6}, where we have that $ y_2 $ is riskier than $ y_1 $ according to the mean-variance (M-V) criterion since $ Var(y_2) > Var(y_1). $ The investment problem described in Table \ref{tab6} has been chosen to illustrate the role of fractional entropy based decision models, viz., EU-FE and EU-FEV models in analyzing decision problems where the Shannon entropy values are equal. Using the current models, the decision choices of this problem are analyzed for risk-neutral utility function $ u(x)= x $ and risk-averse utility functions $ u(x) = \log (x) $ and $ u(x)=\sqrt{x} $ . \\
\hspace*{0.2in} We start by computing the risk measure functions in terms of the risk-tradeoff factor $ \lambda $ for each of the utility functions $u(x)$ using  fractional parameter $ \alpha=0.4 $. The results can be found in Table \ref{tab7}. We get different entropy values corresponding to our models due to the presence of the fractional parameter $ \alpha $. Therefore, from Table \ref{tab7}. we find that the EU-FE and EU-FEV models assign different risk measure values $ 1.70 \lambda - 1 $ and $ 2.21 \lambda - 1 $ to $ y_1 $ and $ y_2, $ respectively, for a risk-neutral utility function $ u(x) = x $. However,  the EU-E model gives $ R(y_1)=R(y_2)=2.47 \lambda - 1 ~ \forall~ \lambda \in [0,1] $. So, an individual cannot reach to a conclusion based on the previous Eu-E model. 
Therefore, for risk-neutral utility function, the individual decisions based on their risk preferences can be analyzed by the fractional entropy function alone without the aid of variance. 
This highlights the advantages of using the fractional entropy function over Shannon entropy function in human decision analysis models. Moreover, from Table \ref{tab8}, we find that the range of $ \lambda $ values supporting a risky decision ($R(y_1) < R(y_2)$) defined on the basis of the M-V criterion  aligns well with the predictions of the fractional entropy-based risk models for all the types of utility functions considered in our study.  For the risk-neutral utility function $ u(x)=x,$ we get $ R(y_1)=R(y_2) $ for $ \lambda = 0 $. In this case, the decision maker shows no specific preference for any action. Further, the variations in the M-V defined decision outcomes ($R(y_1) < R(y_2)$)  are shown with respect to variations in the risk sensitivity parameter $ \alpha $, keeping $ \lambda $ fixed, as illustrated in Fig. \ref{Fig9}.
\begin{r1}
For a decision analysis model $ G(\Theta,Y,u) $ with risk-neutral utility function, $ u(x) = bx + c, ~b>0, $ if $ E(y_1) = E(y_2)~ \text{and} ~ Var(y_1) = Var(y_2) ~ \forall~ y_1,y_2 \in Y, $ then the EU-FE and EU-FEV risk measures can be ordered according to the entropy values, i.e., for a particular value of $ \alpha, $ we have
$$R(y_1) < R(y_2) \iff H_{y_1}^\alpha(\theta) < H_{y_1}^\alpha(\theta)$$  
\end{r1}
\subsubsection{Levy paradoxical problem} 
The decision paradoxical problem, as proposed by Levy (1992) is described in Table \ref{tab9}. The decision results of this problem based on the M-V criterion considers $ y_2 $ to be riskier than $ y_1 $ since $ E[y_1] > E[y_2] $ and $Var[ y_1 ] < Var[ y_2 ].$ Now, we interpret this result in terms of the proposed decision models. For this we have evaluated the risk measure expressions in terms of $ \lambda $ associated with each action for $\alpha = 0.4 $. The computed results are provided in Table \ref{tab10}. From the risk expressions, we obtained the range of $ \lambda $ satisfying the M-V based decision outcome ($ R(y_1)< R(y_2) $). The outputs of $\lambda$ connected to the decision outcomes based on the proposed risk measures are provided in Table \ref{tab11}. \\
\hspace*{0.2in} From the results of Table \ref{tab11}, we find that for a risk-neutral function $ u(x)=x, $ the EU-EV and EU-FEV models gives the same decision outcome ($R( y_1 ) < R( y_2)$ as predicted by the M-V criterion for all values of $ \lambda \in [0,1]. $  However, the EU-FE decision model  prefers the action $ y_1 $ over $ y_2 $ for $ \lambda \in [0,0.07) $. This result is almost similar to that of the EU-E model. However, the decision outcomes based on EU-FE and EU-FEV risk measures or the range of $ \lambda $ for which the action $y_1$ will be preferred over $y_2$ may change with values of $\alpha$ other than 0.4. These adjustments were not possible with the previous models using Shannon entropy, thereby failing to take the individual risk preferences due to uncertainty of outcomes into account. This highlights the superiority of the EU-FEV and EU-FE models, attributed to the inclusion of the sensitivity parameter $ \alpha $ in the Ubriaco entropy function. \\
\hspace*{0.2in}  Further, for the risk averse function $ u(x)=\log(x), $ the EU-EV and EU-FEV risk models give that for $ \alpha=0.4, ~y_2$ is riskier than $ y_1 $ for $ \lambda \in (0.75,1] $ and $ \lambda \in (0.78,1], $ respectively. So there is not much disparity in the results for the EU-FEV model with $ \alpha = 0.4 $ and the EU-EV model. However, both EU-E and EU-FE models indicate that $ y_1 $ is riskier than $ y_2 $ for all values of $ \lambda. $ The same result follows with $ u(x)=\sqrt{x} $ for the EU-E and EU-FE models. 
Further, for a fixed value of the risk tradeoff factor $ \lambda = 0.5 $, the variations in the decision outcomes are also demonstrated with respect to changing values of $ \alpha $ through Fig. \ref{Fig10}. 
\begin{p1}
For a decision analysis model $ G(\Theta,Y,u) $ with risk-neutral utility function, $ u(x) = bx + c, ~b>0, $ if \[E(y_1) > E(y_2) ~ \forall~ y_1,y_2 \in Y\] and for a particular value of $ \alpha, $ \[\frac{Var[y_1]-Var[y_2]}{max_{y\in Y}\{Var[y]\}} + H_{y_1}^\alpha(\theta)-H_{y_2}^\alpha(\theta) < 0,\] then $ R(y_1) < R(y_2) ~\forall~ \lambda \in [0,1] $ with EU-FEV decision model. 
\end{p1}
\begin{r1}
In case of the Levy paradoxical problem, \[\frac{Var[y_1]-Var[y_2]}{max_{y\in Y}\{Var[y]\}} = -0.84 ~\text{and}~ H_{y_1}^{0.4}(\theta)-H_{y_2}^{0.4}(\theta) = 0.5.\] Hence, a decision maker will always choose the action $ y_1 $ under EU-FEV decision model.
\end{r1}
\subsubsection{Allais paradoxical problem}
The  well-known economic decision paradox, introduced by Allais (1953), demonstrates that expected utility alone is insufficient to capture individuals' decision-making behavior in situations involving risk and uncertainty. This paradox involves two experiments where individuals are given the choice between two prospects with potential capital gains, viz., $ y_1 $ or $ y_2 $ in the first experiment and $ y_3 $ or $ y_4 $ in the second one as shown in Table \ref{tab12}. Further, it is well established that decision-makers tend to prefer $ y_1 $ over $ y_2 $ and $ y_4 $ over $ y_3. $ So this preference is compared with the proposed EU-FE and EU-FEV models taking three types of utility functions, viz., risk-neutral ($ u(x)=x $), risk-averse ($ u(x)= \sqrt{x} $) and risk-seeking ($ u(x)=x^2 $) to assess how well the current decision models align with existing models in anticipating decision paradoxes. \\
\hspace*{0.2in} The computed risk measure values in terms of $ \lambda $ are displayed in Table \ref{tab13} to get the range of $ \lambda $ that meets the desired conditions obtained from the M-V criterion. From Table \ref{tab14}, it is evident that the EU-FE model with $ \alpha=0.4 $ demonstrates greater alignment in the range of $ \lambda $ values with the preference patterns derived from the M-V criterion, as compared to the EU-E model. Further, with the inclusion of variance in the risk measure, we get improved results for the first experiment of choosing actions between $ y_1 $ and $ y_2 $ for the risk-neutral and risk-averse utility functions. However, the EU-FE model gives better results than the EU-FEV and EU-EV models as par the M-V criterion outcome for the second experiment of choices between actions $ y_3 $ and $ y_4. $ It also gives better result than the EU-FEV model for the first experiment for the risk-seeking utility function. This highlights the importance of both the fractional entropy-based models with and without the inclusion of variance. Hence, the inclusion of variance is not always necessary with the fractional-entropy based decision models for better predictions of choices. This is due to the dependence of individual decisions on both expected utility and subjective uncertainty of outcomes, which can be portrayed well with EU-FE and EU-FEV decision models. This dependence of the decision outcomes on the risk attitudes towards uncertainty of outcomes is depicted in Fig. \ref{Fig11} which showed the variations of individual preferences with respect to $ \alpha, $ keeping $ \lambda $ fixed. This underscores the importance of using Ubriaco entropy in risk measures for predicting decisions based on the subjective uncertainty of outcomes or states and their subjective utilities, even when the cardinality of the set of outcomes or states differs across actions, as demonstrated in the Allais problem.\\
\hspace*{0.2in} These decision results indicate that the proposed model is highly compatible and outperforms previous Shannon entropy-based EU-E and EU-EV models in many situations for $ \alpha = 0.4 $. This is possible because of the flexibility of the Ubriaco entropy function to accommodate the risk tolerance levels of individuals associated with uncertainty of outcomes to derive desired results by tuning the risk sensitivity parameter $ \alpha $.

\section{Portfolio stock selection using EU-FE and EU-FEV models}
In the previous subsection, the potentials of the proposed fractional Ubriaco entropy-based decision analysis models in solving some important investment problems and decision paradoxes were explored. The decision results revealed that it outperformed the traditional Shannon entropy based decision models in some cases.	
Therefore, the developed EU-FE and EU-FEV models introduced in the previous section have been adopted to select the optimal combination of stock components from the Portuguese Stock Index (PSI 20) index downloaded from  \textit{Investing.com} for generating efficient portfolios. The set comprising of all the efficient portfolios with a given number of stocks gives the efficient frontier. In the current study, the mean-variance portfolio optimization method is adopted after selecting an optimal combination of stocks for generating efficient frontiers as done in Brito (2023), by applying machine learning methods such as Artificial Neural Network (ANN) approach. We have used scaled conjugate gradient algorithm as the training function and the data sets are divided into 70\% training, 15\% validation and 15\% testing. Once trained, the ANN predicts the risk for each stock for all combinations of $ \alpha $ and $ \lambda $ using the formulated fractional order entropy-based risk measures.
The objective of the ANN training is to minimize the mean square of the residuals (MSE), computed as:
\begin{equation}
	MSE = \frac{1}{N} \sum_{j=1}^{N} ({Res}_i)^2,  
\end{equation} where $ N $ is the number of samples and $ {Res}_i = y_i - \hat{y_i};~ y_i = $  actual value of the target variable (the true risk value from our decision model) and $ \hat{y_i} = $ predicted value from the ANN for the same input vector $ x_i $. \\
\hspace*{0.2in} For applying the ANN approach in the stock selection decision modeling, we consider an investor willing to choose a set of $ m $ stock components from $ n $ options. Then the action of selecting stock $ M_i $ will be represented by $ k_i $ for $ i = 1,2,\dots,m $ and the action space will be given by $ K = \{k_1, k_2, \dots, k_n\}. $ Next, the historical data of closing prices for each stock component is collected for PSI 20 data series from 01/01/2019 to 31/12/2020. The stock components included in the PSI 20 data collected are listed in Table \ref{Components}. We obtain a set of 512 entries for each stock $ M_i $ denoted by \{$ p_{i0},p_{i2},\dots,p_{iL} $\}, where $ L = 511. $ Then the returns for each stock $ M_i $ can be obtained as:
\begin{equation}
r_{il} = \log \left(\frac{p_{il}}{p_{i(l-1)}}\right),
\label{4.9}
\end{equation} where i = 1, \dots, 15 and $ l = 1,2,\dots,511. $ Then we obtain a set of returns given by \{$ r_{i1},r_{i2},\dots,r_{i511}, $\} for stock $ M_i $. Among all the stock returns obtained, we find that the minimum return value $ r_{min}= -0.0929 $ and the maximum value $ r_{max}= 0.0751. $ This gives us that 
\[r_{il}\in [-0.0929,0.0751],\] for $ i = 1,2,\dots,15 $ and $ l = 1, \dots, 511. $ The interval [$ r_{min}, r_{max} $] is partitioned into $ J = 15 $ sub-intervals of equal length $ \Delta = \frac{r_{max}-r_{min}}{J} = 0.0112 $. Thus we get 15 sub-intervals of the form $ J_1 = [-0.0929,-0.0817), J_2 = [-0.0817,-0.0705), \dots ,J_{15} = [0.0751,0.0639] $, each of length $ \Delta=0.0112. $ The utility function considered in this study is a S-shaped utility function, known to be concave for profits and convex for losses, as introduced by Kahneman and Tversky (1979) to represent the degree of risk aversion among decision-makers. The utility function is defined as:
\begin{equation}
u(x) = \begin{cases}
	\log (1+x) & x \geq 0,\\
	-\log (1-x) & x < 0.
\end{cases}
\end{equation} Therefore, the EU-FE and EU-FEV risk measures are defined using this utility function for a each stock components. The values of the risk factors used in defining the fractional order entropy based risk measures, given by (\ref{3.8}) and (\ref{3.9}), are provided in Table \ref{RiskValues}. \\
\hspace*{0.2in} Firstly, the set of 15 stocks of the chosen financial data provided in Table \ref{Components} are ranked by implementing a bootstrapped ANN approach on the basis of their predicted risk values based on the proposed decision models (EU-FE and EU-FEV). Higher rank is associated with a stock with lower risk. 100 bootstrap iterations are implemented for training the data with ANN. This creates 100 different predictions for each sample point, allowing us to compute confidence intervals (CIs) over risk predictions. Thus an average risk value is obtained for each stock over 100 bootstrap models. The stocks are then ranked in terms of ascending order of their fractional entropy based risk scores. The rank of the stocks and their corresponding CIs of predicted risk are provided in Table \ref{Ranks}. Stocks with lower predicted risk and narrower CI are better. \\
\hspace*{0.2in} From Fig. \ref{Fig20}, we can observe that the efficient frontier of selected subset of top 7 stocks almost intersects the frontier of portfolios having all 15 stock components of PSI 20 data for return values above 0.5$\times 10^{-3}$. This shows that the efficiency of a portfolio with 7 stocks is almost equal to that of a portfolio set of all stocks for a fixed return range. Thus the EU-FE and EU-FEV models can be effectively used as a stock selection model before the construction of efficient frontiers of portfolios. From this result, we can hereby state that portfolios with less than half the total number of stocks (7 out of 15) in a real stock market data (Portuguese stock market data in this case) can generate equally efficient frontiers as that of portfolios with all the stock components. Thus the current fractional entropy based models can be effectively applied to decrease the number of stocks without affecting the diversification or efficiency of the constructed portfolio.

\section{Conclusion}
The current study highlights that the risk-seeking behavior or tolerance levels of decision makers significantly influence decisions made under uncertainty and risk. The parameter $ \alpha $ in the fractional order entropy proposed by Ubriaco is interpreted as a measure of an individual's risk tolerance or risk sensitivity to uncertainty. Specifically, lower values of $ \alpha $ (closer to zero) indicate lower risk tolerance (conservative attitude), while values approaching one correspond to individuals with higher tolerance for risk (adventurous attitude).
From the findings of this study, it can be remarked that the proposed fractional entropy-based (EU-FE, EU-FEV) decision analysis models is a much better choice for analyzing some important investment problems and economic decision paradoxes. It also serves as an efficient tool for modeling stock selection decisions under risk. The selected stocks can be effectively used to generate efficient frontiers of portfolios using mean-variance optimization method.

\newpage
\section*{Figures and Tables}
\begin{table}[h]
\centering
\tiny
\caption{Portfolio model 1 (Diversified)}
\begin{tabular}{| l | l | c |}
\hline
\textbf{Sl. No.} & ~~~~\textbf{Company (Diversified)} & ~~~~\textbf{Probabilities ($p(x)$)}\\
\hline \hline
1. & ~~~~~~~Pl Industries & ~~~~~0.03 \\
2. & ~~~~~~~Eicher Motors & ~~~~~0.025 \\
3. & ~~~~~~~Balkrishna Industries & ~~~~~0.03 \\
4. & ~~~~~~~Minda Corp & ~~~~~0.03 \\
5. & ~~~~~~~HDFC Bank & ~~~~~0.05 \\
6. & ~~~~~~~HDFC Limited & ~~~~~0.04 \\
7. & ~~~~~~~Bajaj Finance & ~~~~~0.04 \\
8. & ~~~~~~~State Bank of India & ~~~~~0.04 \\
9. & ~~~~~~~Bajaj Finserv & ~~~~~0.04 \\
10. & ~~~~~~~Multi Commodity Exchange & ~~~~~0.03 \\
11. & ~~~~~~~Larsen \& Toubro & ~~~~~0.03 \\
12. & ~~~~~~~KSB Pumps & ~~~~~0.03 \\
13. & ~~~~~~~Thermax & ~~~~~0.03 \\
14. & ~~~~~~~Bharat Electronics & ~~~~~0.04 \\
15. & ~~~~~~~UltraTech Cement & ~~~~~0.2 \\
16. & ~~~~~~~KNR Construction & ~~~~~0.03 \\
17. & ~~~~~~~Jubilant Food & ~~~~~0.025 \\
18. & ~~~~~~~Tata Consumer & ~~~~~0.025 \\
19. & ~~~~~~~Titan Co & ~~~~~0.025 \\
20. & ~~~~~~~Pidilite Industries & ~~~~~0.03 \\
21. & ~~~~~~~Radico Khaitan & ~~~~~0.03 \\
22. & ~~~~~~~Bata India & ~~~~~0.03 \\
23. & ~~~~~~~Infosys & ~~~~~0.04 \\
24. & ~~~~~~~Info Edge & ~~~~~0.02 \\
25. & ~~~~~~~Tata Consultancy Services & ~~~~~0.03 \\
26. & ~~~~~~~Tech Mahindra Limited & ~~~~~0.02 \\
27. & ~~~~~~~TeamLease Services & ~~~~~0.03 \\
28. & ~~~~~~~Container Corporation of India & ~~~~~0.03 \\
29. & ~~~~~~~Reliance Industries & ~~~~~0.04 \\
30. & ~~~~~~~Divis Laboratories & ~~~~~0.025 \\
31. & ~~~~~~~Syngene International & ~~~~~0.03 \\
32. & ~~~~~~~Brigade Enterprises & ~~~~~0.03 \\		
\hline \hline
Total &~~~~~~~ & ~~~~~1 \\
\hline			
\end{tabular}	
\label{tab:di}
\end{table}
\begin{table}[]
\tiny
\centering
\caption{Portfolio model 2 (Large Cap)}
\begin{tabular}{| c | l | c |}
\hline
\textbf{Sl. No.} & ~~~~\textbf{Company (Large Cap)} & ~~~~\textbf{Probabilities ($p(x)$)}\\
\hline \hline
1. & ~~~~~~~HDFC Bank & ~~~~~0.1 \\
2. & ~~~~~~~HDFC Limited & ~~~~~0.08 \\
3. & ~~~~~~~Bajaj Finance & ~~~~~0.08 \\
4. & ~~~~~~~State Bank of India & ~~~~~0.08 \\
5. & ~~~~~~~Eicher Motors & ~~~~~0.05 \\ 
6. & ~~~~~~~Larsen \& Toubro & ~~~~~0.06 \\
7. & ~~~~~~~UltraTech Cement & ~~~~~0.04 \\
8. & ~~~~~~~Jubilant Food & ~~~~~0.05 \\
9. & ~~~~~~~Tata Consumer & ~~~~~0.05 \\
10. & ~~~~~~~Titan Co & ~~~~~0.05 \\
11. & ~~~~~~~Infosys & ~~~~~0.08 \\
12. & ~~~~~~~Info Edge & ~~~~~0.05 \\
13. & ~~~~~~~Tata Consultancy Services & ~~~~~0.06 \\
14. & ~~~~~~~Tech Mahindra Limited & ~~~~~0.04 \\
15. & ~~~~~~~Reliance Industries & ~~~~~0.08 \\
16. & ~~~~~~~Divis Laboratories & ~~~~~0.05 \\
\hline \hline
Total &~~~~~~~ & ~~~~~1 \\
\hline					
\end{tabular}	
\label{tab:lc}
\end{table}
\begin{table}[]
\tiny
\centering
\caption{Portfolio model 3 (Mid Cap)}
\begin{tabular}{| l | l | c |}
\hline
\textbf{Sl. No.} & ~~~~\textbf{Company (Mid Cap)} & ~~~~\textbf{Probabilities ($p(x)$)}\\
\hline \hline
1. & ~~~~~~~Balkrishna Industries & ~~~~~0.06 \\
2. & ~~~~~~~Pl Industries & ~~~~~0.06 \\
3. & ~~~~~~~Minda Corp & ~~~~~0.06 \\
4. & ~~~~~~~Bajaj Finserv & ~~~~~0.08 \\
5. & ~~~~~~~Multi Commodity Exchange & ~~~~~0.06 \\
6. & ~~~~~~~KSB Pumps & ~~~~~0.06 \\
7. & ~~~~~~~Thermax & ~~~~~0.06 \\
8. & ~~~~~~~Bharat Electronics & ~~~~~0.08 \\
9. & ~~~~~~~Pidilite Industries & ~~~~~0.06 \\
10. & ~~~~~~~Radico Khaitan & ~~~~~0.06 \\
11. & ~~~~~~~Bata India & ~~~~~0.06 \\
12. & ~~~~~~~TeamLease Services & ~~~~~0.06 \\
13. & ~~~~~~~Container Corporation of India & ~~~~~0.06 \\
14. & ~~~~~~~Syngene International & ~~~~~0.06 \\
15. & ~~~~~~~Brigade Enterprises & ~~~~~0.06 \\
16. & ~~~~~~~KNR Construction & ~~~~~0.06 \\
\hline \hline
Total &~~~~~~~ &  ~~~~~1 \\
\hline					
\end{tabular}	
\label{tab:mc}
\end{table}	
\begin{table}[]
\tiny
\centering
\caption{Portfolio model 4 (Hypothetical)}
\begin{tabular}{| l | l | c |}		
\hline
\textbf{Sl. No.} & ~~~~\textbf{Company} & ~~~~\textbf{Probabilities ($p(x)$)}\\
\hline \hline
1. & ~~~~~~~Guj. State Petro Ltd. & ~~~~~0.5 \\
2. & ~~~~~~~Birla Corporation & ~~~~~0.04 \\
3. & ~~~~~~~Oberoi Realty & ~~~~~0.2 \\
4. & ~~~~~~~Tata Chemicals & ~~~~~0.06 \\
5. & ~~~~~~~Century Textiles & ~~~~~0.2 \\
\hline \hline
Total &~~~~~~~ & ~~~~~1 \\
\hline						
\end{tabular}	
\label{tab:hy}
\end{table}
\begin{table}[]
\small
\centering
\caption{Values of Ubriaco entropy for different portfolios}
\begin{tabular}{| l | c  c  c  c  c  c  |}
\hline
\textbf{Type} & $ \alpha = 0.1 $ ~& $ \alpha = 0.3 $ ~& $ \alpha = 0.5 $ ~& $ \alpha = 0.7 $ ~& $ \alpha = 0.9 $  ~& Shannon  \\
\hline \hline
Diversified (Di) & 1.1314 & 1.4486 & 1.8549 & 2.3756 & 3.0429 & 3.4440 \\
\hline	
Large-cap (Lc) & 1.1053 & 1.3508 & 1.6515 & 2.0199 & 2.4716 & 2.7344 \\
\hline
Mid-cap (Mc) & 1.1071 & 1.3569 & 1.6632 & 2.0388 & 2.4994 & 2.7674 \\
\hline
Hypothetical (Hy) & 0.9319 & 0.816 & 0.7224 & 0.6467 & 0.5854 & 0.5593 \\
\hline			
\end{tabular}
\label{tab:entropy}	
\end{table}
\begin{table}[h]
\small
\centering
\caption{Ordering of entropy values for different portfolios}
\begin{tabular}{| l | l |}
\hline
Di portfolio:&~ $ H_U^{0.1} < H_U^{0.3} < H_U^{0.5} < H_U^{0.7} < H_U^{0.9} < $ Shannon \\
Lc portfolio:&~ $ H_U^{0.1} < H_U^{0.3} < H_U^{0.5} < H_U^{0.7} < H_U^{0.9} < $ Shannon  \\
Mc portfolio:&~ $ H_U^{0.1} < H_U^{0.3} < H_U^{0.5} < H_U^{0.7} < H_U^{0.9} < $ Shannon  \\
Hy portfolio:&~ Shannon $ < H_U^{0.9} < H_U^{0.7} < H_U^{0.5} < H_U^{0.3} < H_U^{0.1} $.\\	
\hline
\end{tabular}
\label{tab:order1} 	
\end{table}    

\begin{table}[h]
\tiny
\centering
\caption{Nawrocki and Harding problem with two different actions}
\vspace*{0.1in}
\begin{tabular}{ c c c c c c c c c c }
\hline
&~  &~  $\theta_{1}$ &~ $\theta_{2}$ &~ $\theta_{3}$ &~ $\theta_{4}$ &~ $\theta_{5}$ &~ Expected value &~ Variance &~ Entropy \\
\hline \hline 
$ y_1 $ &~ $x_{1j}$ &~ 1 &~ 2 &~ 3 &~ 4 &~ 5 &~ 3 &~ 1.2 &~ 1.47 \\
&~ $ p_{1j} $ &~ 0.1 &~ 0.2 &~ 0.4 &~ 0.2 &~ 0.1 &~  &~  &~ \\
\hline
\vspace*{0.05in}
$ y_2 $ &~ $ x_{2j} $ &~ 1 &~ 2 &~ 3 &~ 4 &~ 5 &~ 3 &~ 1.8 &~ 1.47 \\
&~ $ p_{1j} $ &~ 0.2 &~ 0.1 &~ 0.4 &~ 0.1 &~ 0.2 &~  &~  &~ \\
\hline \hline			 
\end{tabular}
\label{tab6}
\end{table}
\begin{table}[h]
\tiny
\centering
\caption{Risk measures for Nawrocki and Harding problem when $ \alpha = 0.4 $}
\vspace*{0.1in}
\begin{tabular}{ c c c c c c c}
\hline
Utility function &~ $ H^{0.4}(\theta) $ &~ $R(y), ~y = a ~ \text{or} ~ b $ &~ EU-E &~ EU-EV &~ EU-FE &~ EU-FEV \\
\hline \hline
\vspace*{0.05in}
$ u(x) = x $ &~ 0.70 &~ $R(y_1)$ &~ $2.47 \lambda -1 $ &~ $2.07\lambda -1 $ &~ $1.70 \lambda - 1 $	&~ $1.68 \lambda - 1 $ \\
&~ 1.21 &~ $R(y_2)$ &~ $2.47 \lambda -1 $ &~ $2.24\lambda -1 $ &~ $2.21 \lambda - 1 $	&~ $1.94 \lambda - 1 $ \\
\hline
\vspace*{0.05in}
$ u(x) = \log (x) $ &~ 0.70 &~ $R(y_1)$ &~ $2.47 \lambda -1 $ &~ $2.07\lambda -1 $ &~ $1.70 \lambda - 1 $ &~ $1.68 \lambda - 1 $ \\
&~ 1.21 &~ $R(y_2)$ &~ $2.42 \lambda -0.95 $ &~ $2.19\lambda -0.95 $ &~ $2.16 \lambda - 0.95 $	&~ $2.06 \lambda - 0.95 $ \\
\hline
\vspace*{0.05in}
$ u(x) = \sqrt{x} $ &~ 0.70 &~ $R(y_1)$ &~ $2.47 \lambda -1 $ &~ $2.07\lambda -1 $ &~ $1.70 \lambda - 1 $ &~ $1.68 \lambda - 1 $ \\
&~ 1.21 &~ $R(y_2)$ &~ $2.46 \lambda -0.99 $ &~ $2.23\lambda -0.99 $ &~ $2.20 \lambda - 0.99 $	&~ $2.10 \lambda - 0.99 $ \\
\hline \hline		 
\end{tabular}
\label{tab7}
\end{table}
\begin{table}[h]
\tiny
\centering
\caption{Risk tradeoff factors for EU-FE and EU-FEV models corresponding to Nawrocki and Harding decision problem}
\vspace*{0.1in}
\begin{tabular}{ c c c c c c }
\hline
Utility function &~ Risk decision &~ EU-E &~ EU-EV &~ EU-FE &~ EU-FEV \\
\hline \hline
\vspace*{0.05in}
$ u(x) = x $ &~ $R(y_1) < R(y_2)$ &~ - &~ $ 0 < \lambda \leq 1 $ &~ $ 0 < \lambda \leq 1 $ &~ $ 0 < \lambda \leq 1 $ \\
&~ $R(y_1) = R(y_2)$ &~ $ 0\leq\lambda\leq 1 $ &~ $\lambda=0 $ &~ $\lambda=0 $ &~ $\lambda=0 $ \\
\hline
\vspace*{0.05in}
$ u(x) = \log (x) $ &~ $R(y_1) < R(y_2)$ &~ $ 0 \leq \lambda < 1$ &~ $ 0 \leq \lambda \leq 1 $ &~ $ 0 \leq \lambda \leq 1 $ &~ $ 0 \leq \lambda \leq 1 $ \\
&~ $R(y_1) = R(y_2)$ &~ $ \lambda = 1 $ &~ - &~ - &~ - \\
\hline
\vspace*{0.05in}
$ u(x) = \sqrt{x} $ &~ $R(y_1)< R(y_2)$ &~ $ 0 \leq \lambda < 1 $ &~ $ 0 \leq \lambda \leq 1 $ &~ $ 0 \leq \lambda \leq 1 $ &~ $ 0 \leq \lambda \leq 1 $\\
&~ $R(y_1) = R(y_2)$ &~ $ \lambda = 1 $ &~ - &~ - &~ - \\		
\hline \hline			 
\end{tabular}
\label{tab8}
\end{table}

\begin{table}[h]
\tiny
\centering
\caption{Levy problem with two actions and two states of nature}
\vspace*{0.1in}
\begin{tabular}{ c c c c c c c  }
\hline
&~  &~  $\theta_{1}$ &~ $\theta_{2}$ &~ Expected value &~ Variance &~ Entropy \\
\hline \hline 
$ y_1 $ &~ $x_{1j}$ &~ 1 &~ 100 &~ 20.8 &~ 1568 &~ 0.5 \\
&~ $ p_{1j} $ &~ 0.8 &~ 0.2 &~ &~ &~ \\
\hline
\vspace*{0.05in}
$ y_2 $ &~ $ x_{2j} $ &~ 10 &~ 1000 &~ 19.9 &~ 9703 &~ 0.06 \\
&~ $ p_{2j} $ &~ 0.99 &~ 0.01 &~ &~ &~  \\
\hline \hline			 
\end{tabular}
\label{tab9}
\end{table}
\begin{table}[h]
\tiny
\centering
\caption{Risk measures for $ \alpha=0.4 $ corresponding to Levy problem}
\vspace*{0.1in}
\begin{tabular}{ c c c c c c c}
\hline
Utility function &~ $ H^{0.4}(\theta) $ &~ $R(y), ~y = a ~ \text{or} ~ b $ &~ EU-E &~ EU-EV &~ EU-FE &~ EU-FEV \\
\hline \hline
\vspace*{0.05in}
$ u(x) = x $ &~ 0.68 &~ $R(y_1)$ &~ $1.50 \lambda -1 $ &~ $1.33\lambda -1 $ &~ $1.68 \lambda - 1 $	&~ $1.42 \lambda - 1 $ \\
&~ 0.18 &~ $R(y_2)$ &~ $1.02 \lambda -0.96 $ &~ $1.49\lambda -0.96 $ &~ $1.14 \lambda - 0.96 $	&~ $1.55 \lambda - 0.96 $ \\
\hline
\vspace*{0.05in}
$ u(x) = \log (x) $ &~ 0.68 &~ $R(y_1)$ &~ $0.89 \lambda -0.39 $ &~ $0.72\lambda -0.39 $ &~ $1.07 \lambda - 0.39 $ &~ $0.81 \lambda - 0.39 $ \\
&~ 0.18 &~ $R(y_2)$ &~ $1.06 \lambda -1 $ &~ $1.53\lambda -1 $ &~ $1.18 \lambda - 1 $	&~ $1.59 \lambda - 1 $ \\
\hline
\vspace*{0.05in}
$ u(x) = \sqrt{x} $ &~ 0.68 &~ $R(y_1)$ &~ $1.31 \lambda -0.81 $ &~ $1.14\lambda -0.81 $ &~ $1.49 \lambda - 0.81 $ &~ $1.23 \lambda - 0.81 $ \\
&~ 0.18 &~ $R(y_2)$ &~ $1.06 \lambda -1 $ &~ $1.53\lambda -1 $ &~ $1.18 \lambda - 1 $ &~ $1.59 \lambda - 1 $ \\
\hline \hline		 
\end{tabular}
\label{tab10}
\end{table}	
\begin{table}[h]
\tiny
\centering
\caption{Risk tradeoff factors for EU-FE and EU-FEV models due to Levy decision problem}
\vspace*{0.1in}
\begin{tabular}{ c c c c c c }
\hline
Utility function &~ Risk decision &~ EU-E &~ EU-EV &~ EU-FE &~ EU-FEV \\
\hline \hline
\vspace*{0.05in}
$ u(x) = x $ &~ $R(y_1) < R(y_2)$ &~ $ 0 \leq \lambda < 0.08 $ &~ $ 0 \leq \lambda \leq 1 $ &~ $ 0 \leq \lambda < 0.07 $ &~ $ 0 \leq \lambda \leq 1 $ \\
&~ $R(y_1) > R(y_2)$ &~ $ 0.08 < \lambda\leq 1 $ &~ - &~ $ 0.07 < \lambda \leq 1 $ &~ -\\
\hline
\vspace*{0.05in}
$ u(x) = \log (x) $ &~ $R(y_1) < R(y_2)$ &~ - &~ $ 0.75 < \lambda \leq 1 $ &~ - &~ $ 0.78 < \lambda \leq 1 $ \\
&~ $R(y_1) > R(y_2)$ &~ $ 0 \leq \lambda \leq 1$ &~ $ 0 \leq \lambda < 0.75 $ &~ $ 0 \leq \lambda \leq 1 $ &~ $ 0 \leq \lambda < 0.78 $  \\
\hline
\vspace*{0.05in}
$ u(x) = \sqrt{x} $ &~ $R(y_1)< R(y_2)$ &~ - &~ $ 0.49 < \lambda \leq 1 $ &~ -  &~ $ 0.53 < \lambda \leq 1 $ \\
&~ $R(y_1) > R(y_2)$ &~ $ 0 \leq \lambda \leq 1 $ &~ $ 0 \leq \lambda < 0.49 $ &~ $ 0 \leq \lambda \leq 1 $ &~ $ 0 \leq \lambda < 0.53 $ \\		
\hline \hline			 
\end{tabular}
\label{tab11}
\end{table}

\begin{table}[h]
\tiny
\centering
\caption{Allais problem with four actions and three states}
\vspace*{0.1in}
\begin{tabular}{ c c c c c c c c }
\hline
&~  &~  $\theta_{1}$ &~ $\theta_{2}$ &~ $\theta_{3}$ &~ Expected value &~ Variance &~ Entropy \\
\hline \hline 
$ y_1 $ &~ $x_{1j}$ &~ 1 &~  &~  &~  &~  &~  \\
&~ $ p_{1j} $ &~ 1 &~  &~  &~ 1 &~ 0 &~ 0\\
\hline
\vspace*{0.05in}
$ y_2 $ &~ $ x_{2j} $ &~ 1 &~ 5 &~ 0 &~  &~  &~ \\
&~ $ p_{2j} $ &~ 0.89 &~ 0.1 &~ 0.01 &~ 1.39 &~ 1.46 &~ 0.38  \\
\hline
\vspace*{0.05in}
$ y_3 $ &~ $ x_{3j} $ &~ 1 &~  &~ 0 &~  &~  &~  \\
&~ $ p_{2j} $ &~ 0.11 &~  &~ 0.89 &~ 0.11 &~ 0.1 &~ 0.35  \\
\hline
\vspace*{0.05in}
$ y_4 $ &~ $ x_{4j} $ &~  &~ 5 &~ 0 &~  &~  &~ \\
&~ $ p_{2j} $ &~ &~ 0.1 &~ 0.9 &~ 0.5 &~ 2.25 &~ 0.33  \\
\hline \hline			 
\end{tabular}
\label{tab12}
\end{table}
\begin{table}[h]
\tiny
\centering
\caption{Risk measures for $ \alpha=0.4 $ corresponding to Allais problem}
\vspace*{0.1in}
\begin{tabular}{ c c c c c c c}
\hline
Utility function &~ $ H^{0.4}(\theta) $ &~ $R(y), ~y = a ~ \text{or} ~ b $ &~ EU-E &~ EU-EV &~ EU-FE &~ EU-FEV \\
\hline \hline
&~ &~ &~ &~ &~ &~ \\
$ u(x) = x $ &~ 0 &~ $R(y_1)$ &~ $0.72 \lambda - 0.72 $ &~ $0.72 \lambda - 0.72 $ &~ $0.72 \lambda - 0.72 $ &~ $0.72 \lambda - 0.72 $ \\
&~ 0.54 &~ $R(y_2)$ &~ $1.38 \lambda -1 $ &~ $1.69\lambda -1 $ &~ $1.54 \lambda - 1 $	&~ $1.59 \lambda - 1 $ \\
&~ 0.53 &~ $R(y_3)$ &~ $0.43 \lambda -0.08 $ &~ $0.27\lambda -0.08 $ &~ $0.61 \lambda - 0.08 $	&~ $0.36 \lambda - 0.08 $ \\
&~ 0.51 &~ $R(y_4)$ &~ $0.69 \lambda -0.36 $ &~ $1.02\lambda -0.36 $ &~ $0.86 \lambda - 0.36 $	&~ $1.11 \lambda - 0.36 $ \\
&~ &~ &~ &~ &~ &~ \\
\hline
&~ &~ &~ &~ &~ &~ \\
$ u(x) = \sqrt{x} $ &~ 0 &~ $R(y_1)$ &~ $0.9 \lambda - 0.9 $ &~ $ 0.9 \lambda - 0.9 $ &~ $0.9 \lambda - 0.9 $ &~ $0.9 \lambda - 0.9 $\\
&~ 0.54 &~ $R(y_2)$ &~ $1.38 \lambda -1 $ &~ $1.69\lambda -1 $ &~ $1.54 \lambda - 1 $ &~ $1.59 \lambda - 1 $ \\
&~ 0.53 &~ $R(y_3)$ &~ $0.44 \lambda -0.1 $ &~ $0.29\lambda -0.1 $ &~ $0.63 \lambda - 0.1 $ &~ $0.39 \lambda - 0.1 $ \\
&~ 0.51 &~ $R(y_4)$ &~ $0.53 \lambda -0.2 $ &~ $0.86\lambda -0.2 $ &~ $0.71 \lambda - 0.2 $ &~ $0.96 \lambda - 0.2 $ \\
&~ &~ &~ &~ &~ &~ \\
\hline
&~ &~ &~ &~ &~ &~ \\
$ u(x) = x^2 $ &~ 0 &~ $R(y_1)$ &~ $0.29 \lambda - 0.29 $ &~ $ 0.29 \lambda - 0.29 $ &~ $0.29 \lambda - 0.29 $ &~ $0.29 \lambda - 0.29 $\\
&~ 0.54 &~ $R(y_2)$ &~ $1.38 \lambda -1 $ &~ $1.69\lambda -1 $ &~ $1.54 \lambda - 1 $ &~ $1.37 \lambda - 1 $ \\
&~ 0.53 &~ $R(y_3)$ &~ $0.38 \lambda -0.03 $ &~ $0.22\lambda -0.03 $ &~ $0.56 \lambda - 0.03 $ &~ $0.32 \lambda - 0.03 $ \\
&~ 0.51 &~ $R(y_4)$ &~ $1.06 \lambda -0.74 $ &~ $1.40\lambda -0.74 $ &~ $1.24 \lambda - 0.74 $ &~ $1.49 \lambda - 0.74 $ \\
&~ &~ &~ &~ &~ &~ \\
\hline \hline		 
\end{tabular}
\label{tab13}
\end{table}	
\begin{table}[h]
\tiny
\centering
\caption{Risk tradeoff factors for EU-FE and EU-FEV models corresponding to Allais problem}
\vspace*{0.1in}
\begin{tabular}{ c c c c c c }
\hline
Utility function &~ Risk decision &~ EU-E &~ EU-EV &~ EU-FE &~ EU-FEV \\
\hline \hline
&~ &~ &~ &~ &~ \\
\vspace*{0.05in}
$ u(x) = x $ &~ $R(y_1) < R(y_2)$ &~ $ 0.42 < \lambda \leq 1 $ &~ $ 0.28 < \lambda \leq 1 $ &~ $ 0.34 < \lambda \leq 1 $ &~ $ 0.32 < \lambda \leq 1 $ \\
&~ $R(y_4) < R(y_3)$ &~ $ 0 \leq \lambda\leq 1 $ &~ $0 \leq \lambda < 0.37 $ &~ $ 0 \leq \lambda \leq 1 $ &~ $ 0\leq \lambda < 0.37 $ \\
\hline
&~ &~ &~ &~ &~ \\
\vspace*{0.05in}
$ u(x) = \sqrt x $ &~ $R(y_1) < R(y_2)$ &~ $ 0.21 < \lambda \leq 1 $ &~ $ 0.13 < \lambda \leq 1 $ &~ $ 0.16 < \lambda \leq 1 $ &~ $ 0.14 < \lambda \leq 1 $ \\
&~ $R(y_4) < R(y_3)$ &~ $ 0 \leq \lambda \leq 1$ &~ $ 0 \leq \lambda < 0.18 $ &~ $ 0 \leq \lambda \leq 1 $ &~ $ 0 \leq \lambda < 0.18 $  \\
\hline
&~ &~ &~ &~ &~ \\
\vspace*{0.05in}
$ u(x) = x^2 $ &~ $R(y_1)< R(y_2)$ &~ $ 0.65 < \lambda \leq 1 $ &~ $ 0.51 < \lambda \leq 1 $ &~ $ 0.57 < \lambda \leq 1 $ &~ $ 0.66 < \lambda \leq 1 $ \\
&~ $R(y_4) > R(y_3)$ &~ $ 0 \leq \lambda \leq 1 $ &~ $ 0 \leq \lambda < 0.60 $ &~ $ 0 \leq \lambda \leq 1 $ &~ $ 0 \leq \lambda < 0.61 $ \\	
&~ &~ &~ &~ &~ \\	
\hline \hline			 
\end{tabular}
\label{tab14}
\end{table}
\begin{table}[h!]
\centering
\caption{Stock components of PSI20}
\begin{tabular}{|c|c|c|c|c|c|}
\hline
Stock & Components & Stock & Components & Stock & Components\\
\hline
$M_1$ & ALSS & $M_6$ & EDPR & $M_{11}$ & NOS \\
$M_2$ & BCP & $M_7$ & GALP & $M_{12}$ & NVGR \\
$M_3$ & CORA & $M_8$ & IBS & $M_{13}$ & RENE \\
$M_4$ & CTT & $M_9$ & JMT & $M_{14}$ & SEM \\
$M_5$ & EDP & $M_{10}$ & MOTA & $M_{15}$ & YSO \\ 
\hline
\end{tabular}
\label{Components}
\end{table}

\begin{table}[h!]
\centering
\caption{Risk factors of EU-FE and EU-FEV decision models}
\begin{tabular}{|c|c|c|r|}
\hline
Stocks $M_i$ & NEU($a_i$) =$\frac{E[u(x_{in})]}{max{E[u(x_{in})]}}$ & NV($a_i$) $=\frac{V(x_{in})}{max{V(x_{in})}}$ & $H^\alpha(X_i)$\\
\hline
$M_1$ & -0.2386  & 0.6712 & 1.189223 \\
$M_2$ & -0.7733 & 0.7184 & 1.138423 \\
$M_3$ & 0.1676 & 0.3313 & 1.218279 \\
$M_4$ & -0.3357 & 0.5658 & 1.233337 \\
$M_5$ & 0.4833 & 0.3486 & 1.097892 \\
$M_6$ & 1.0000 & 0.3499 & 1.161922 \\
$M_7$ & -0.6201 & 0.5955 & 1.086162 \\
$M_8$ & -0.6936 & 1.0000 & 0.983342 \\
$M_9$ & 0.1974 & 0.3002 & 1.137902 \\
$M_{10}$ & -0.4411 & 0.9628 & 1.058261 \\
$M_{11}$ & -0.7326 & 0.3635 & 1.161297 \\
$M_{12}$ & -0.4576 & 0.4529 & 1.225268 \\
$M_{13}$ & -0.0905 & 0.1576 & 1.068735 \\
$M_{14}$ & -0.4427 & 0.4243 & 1.209417 \\
$M_{15}$ & -0.2730 & 0.3809 & 1.237089 \\
\hline
\end{tabular}
\label{RiskValues}
\end{table}

\begin{table}
\centering
\caption{Ranks of stock components of PSI 20 data}
\begin{tabular}{|c|c|c|c|}
\hline
Ranks & Stock & Mean Risk & Confidence Interval (CI) \\
\hline
1 &  $ M_6 $ &  0.1152 &  [0.0203, 0.2083] \\
2 &  $ M_5 $ &  0.3509 &  [0.2793, 0.4298] \\
3 &  $ M_9 $ &  0.5096 &  [0.4516, 0.5695]\\
4 &  $ M_3 $ &  0.5449 &  [0.4791, 0.6006]\\
5 &  $ M_{13} $ &  0.6118 &  [0.5475, 0.6980]\\
6 &  $ M_1 $ &  0.7616 &  [0.6952, 0.8340]\\
7 &  $ M_{10} $ &  0.7751 &  [0.7181, 0.8475]\\
8 &  $ M_{15} $ &  0.8042 &  [0.7370, 0.8589]\\
9 &  $ M_4 $ &  0.8209 &  [0.7697, 0.8740]\\
10 & $ M_8 $ &  0.8578 &  [0.7992, 0.9161]\\
11 &  $ M_{14} $ &  0.8762 &  [0.8270, 0.9221]\\
12 & $ M_{12} $ &  0.8806 &  [0.8354, 0.9279]\\
13 & $ M_{7} $ &  0.9198 &  [0.8670, 0.9772]\\
14 & $ M_{2} $ &  0.9601 &  [0.8930, 1.0177]\\
15 & $ M_{11} $ &  0.9780 &  [0.9140, 1.0399]	\\
\hline
\end{tabular}
\label{Ranks}
\end{table}	
\begin{figure}[h]
\centering
\includegraphics[width=0.65\textwidth]{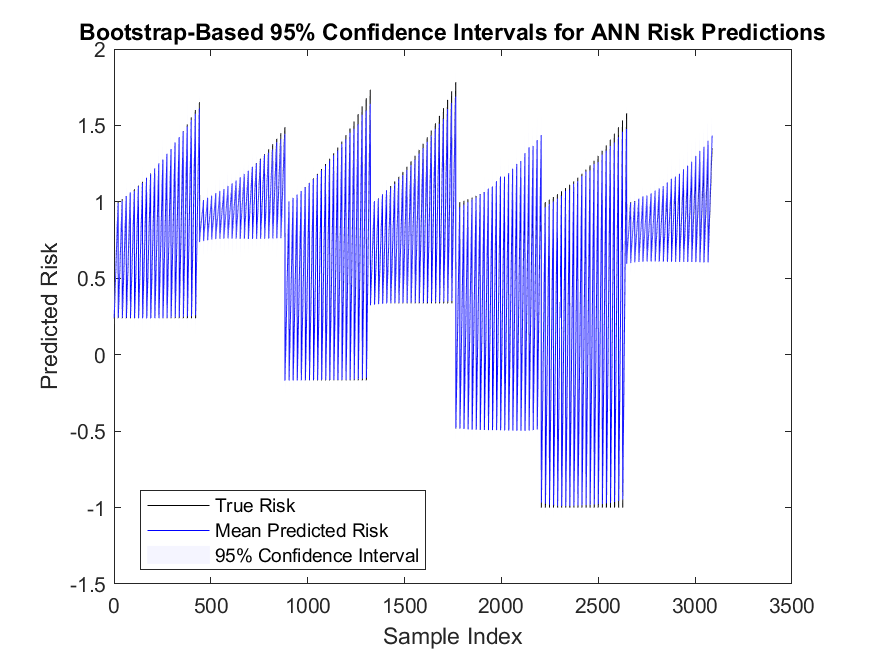}
\caption{Comparison of Efficient Frontiers of portfolio of 7 stocks with that of a portfolio with all the 15 stock components}
\label{Fig21} 
\end{figure} 

\begin{figure}[h]
	\centering
	\includegraphics[width=0.65\textwidth]{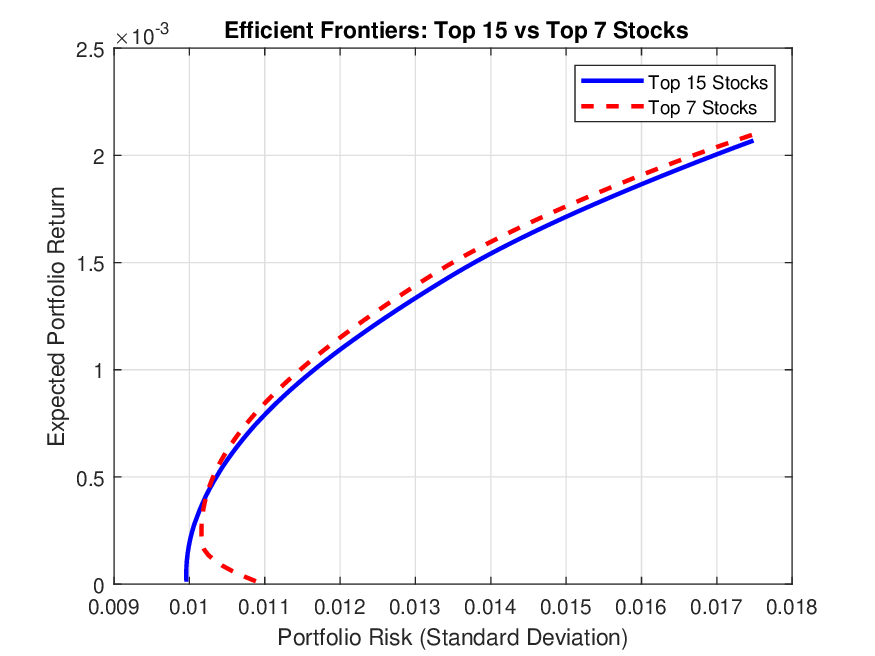} 
	\caption{Comparison of Efficient Frontiers of portfolio of 7 stocks with that of a portfolio with all the stock components}
	\label{Fig20} 
\end{figure}  

\begin{figure}[h]
\centering
\includegraphics[width=0.4\textwidth]{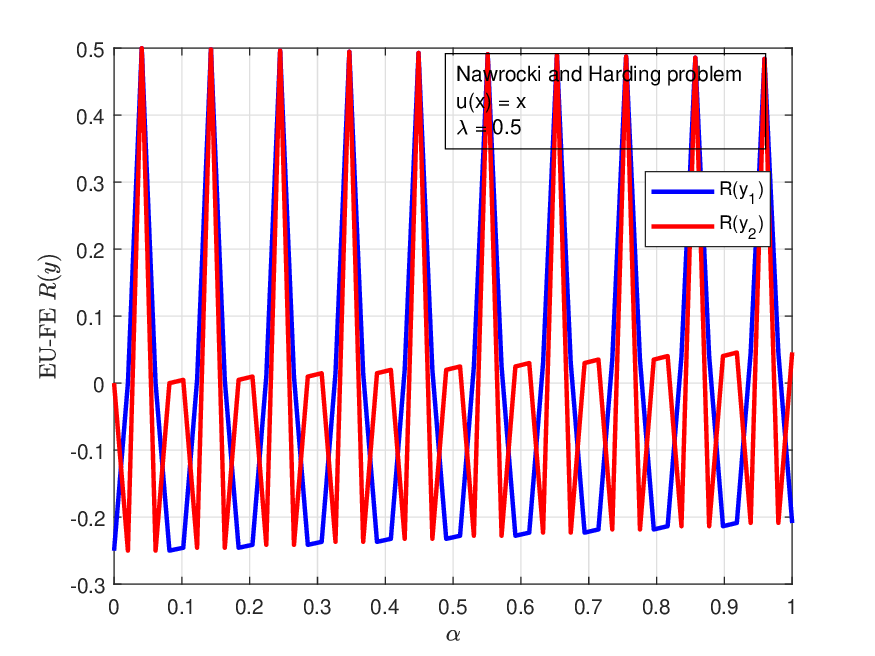} ~~ \includegraphics[width=0.4\textwidth]{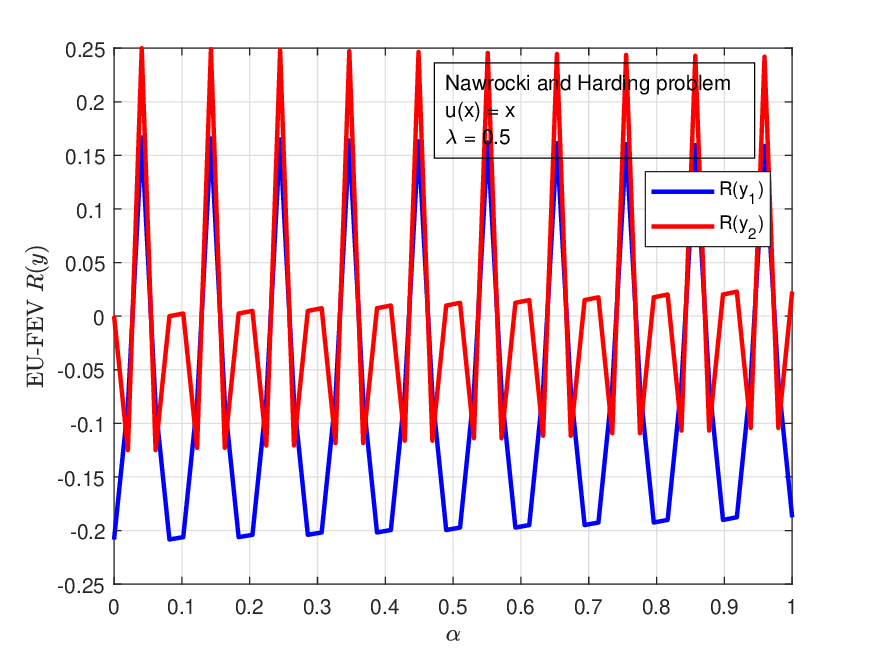} \\
\includegraphics[width=0.4\textwidth]{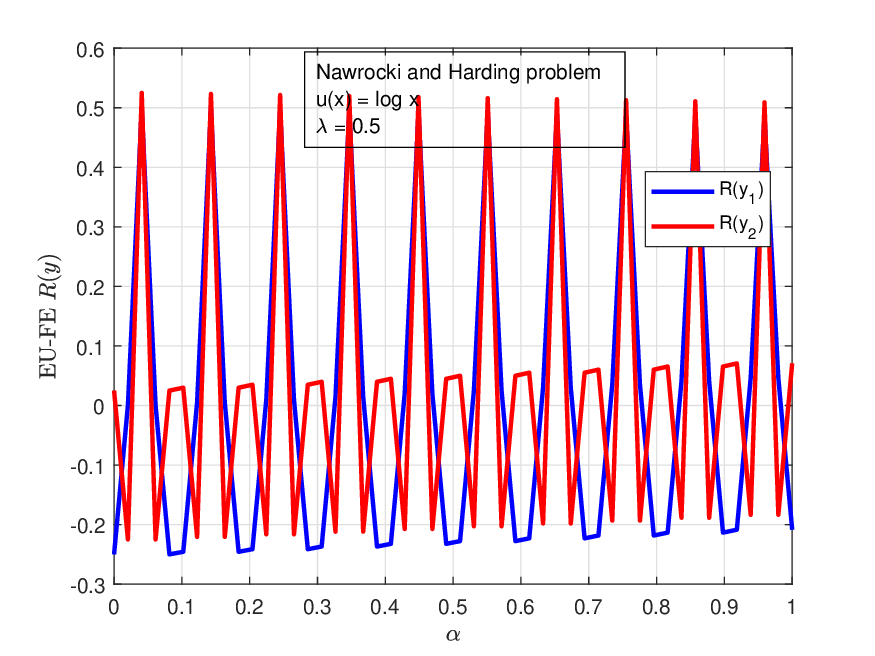} ~~ \includegraphics[width=0.4\textwidth]{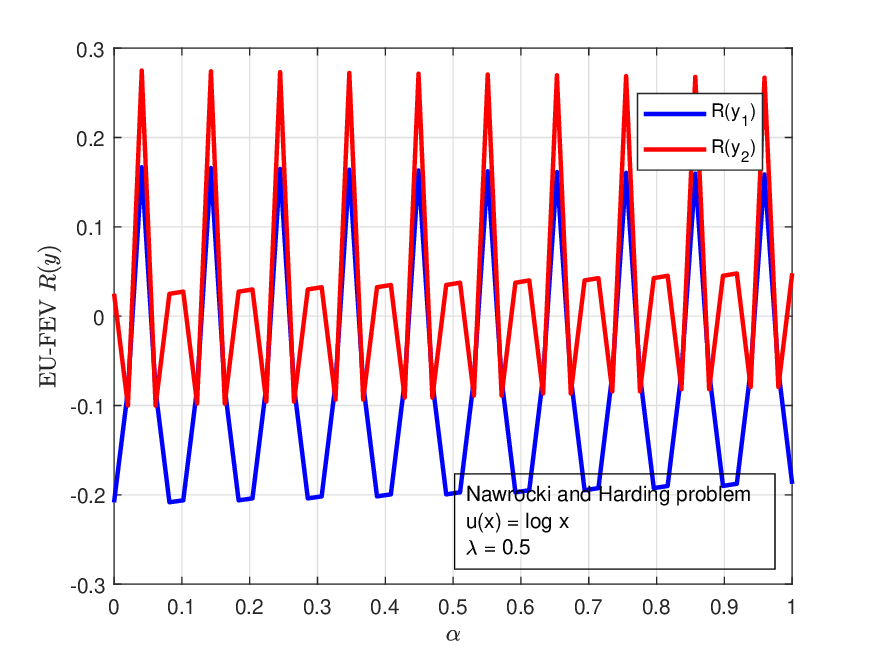} \\
\includegraphics[width=0.4\textwidth]{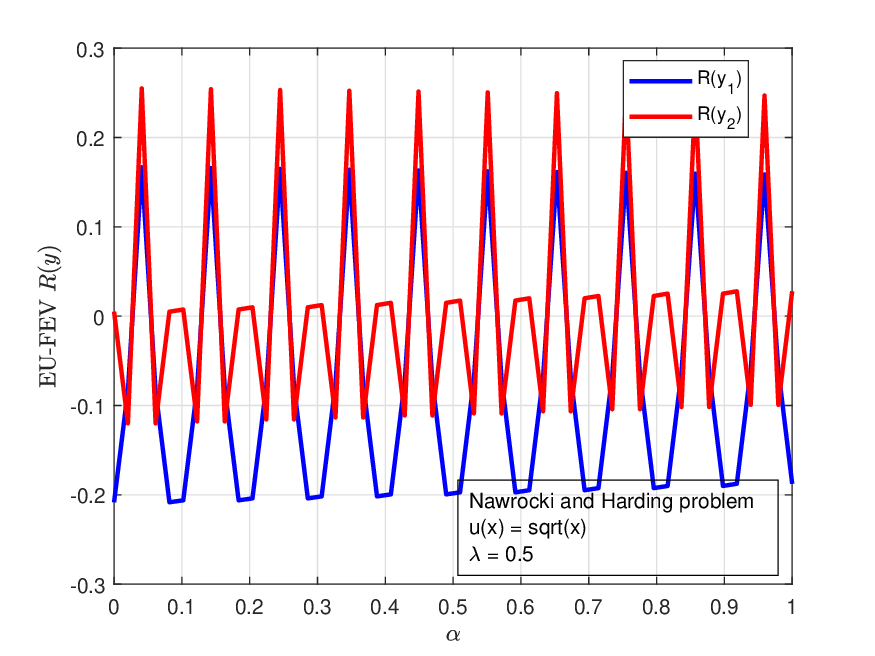} ~~ \includegraphics[width=0.4\textwidth]{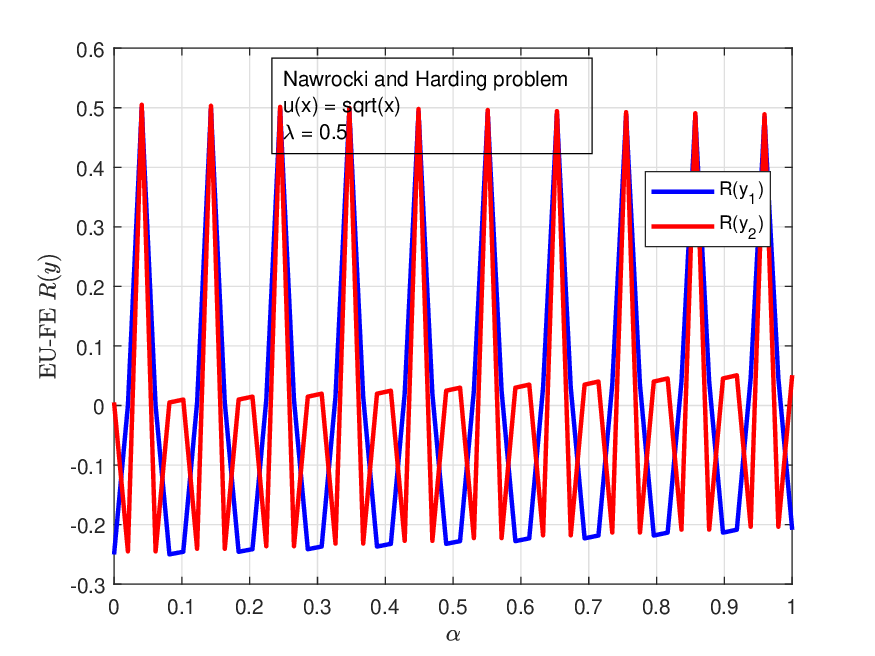}
\caption{Risk measures for actions $ y_3 $ and $ y_4 $ for different $ \alpha $ and $ u(x) $ in Nawrocki and Harding problem}
\label{Fig9} 
\end{figure} 
\begin{figure}[h]
\centering
\includegraphics[width=0.4\textwidth]{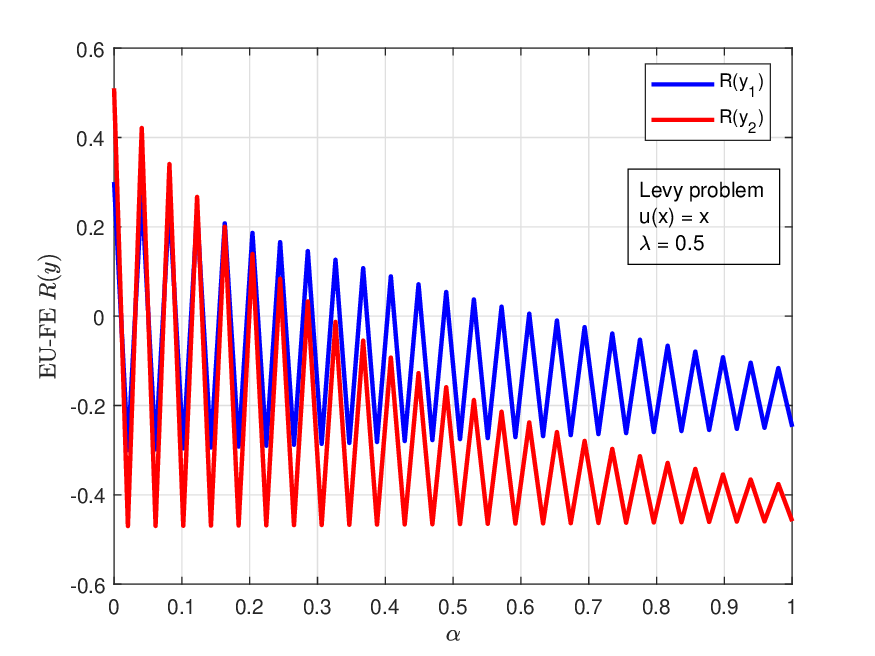} ~~ \includegraphics[width=0.4\textwidth]{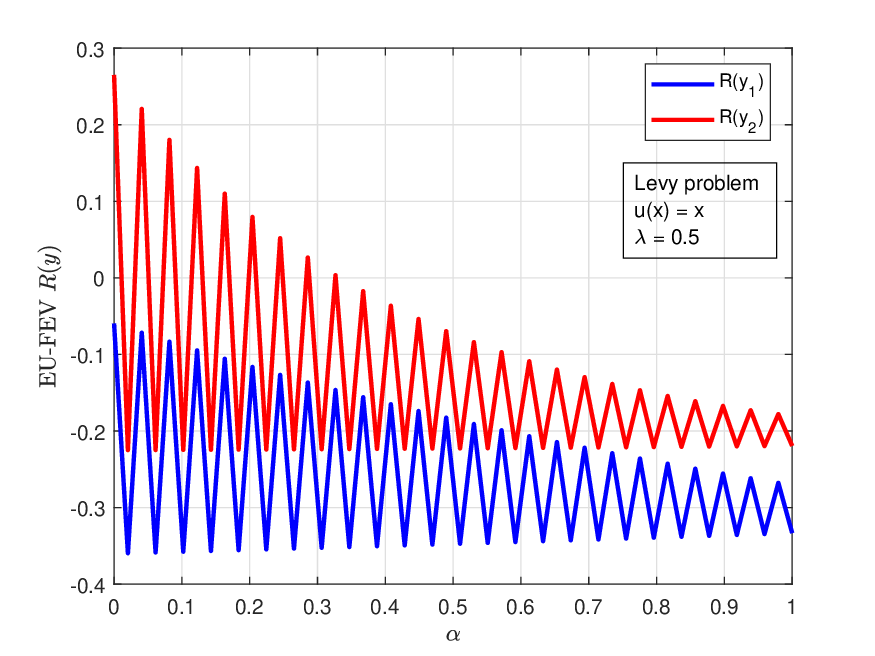} \\
\includegraphics[width=0.4\textwidth]{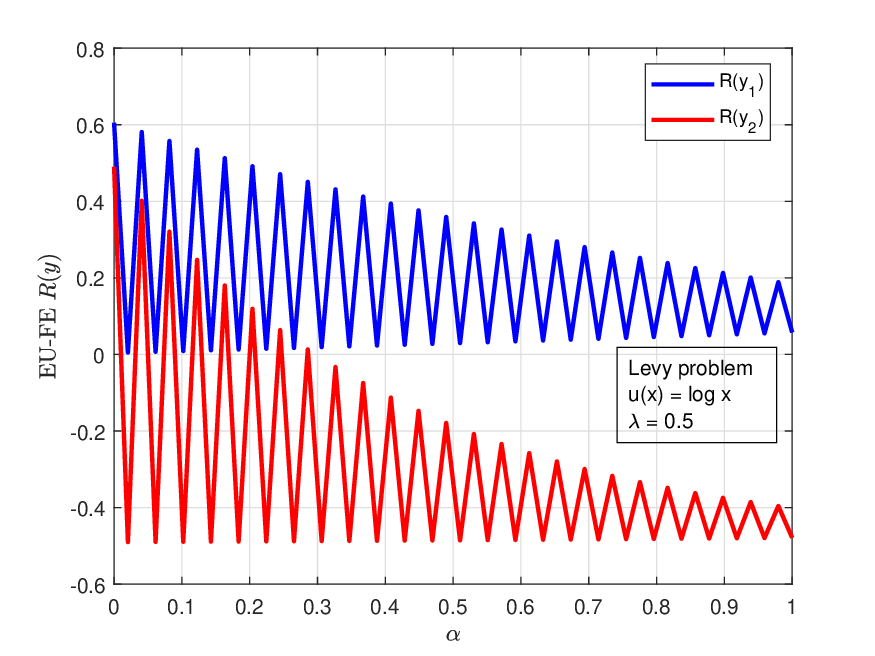} ~~ \includegraphics[width=0.4\textwidth]{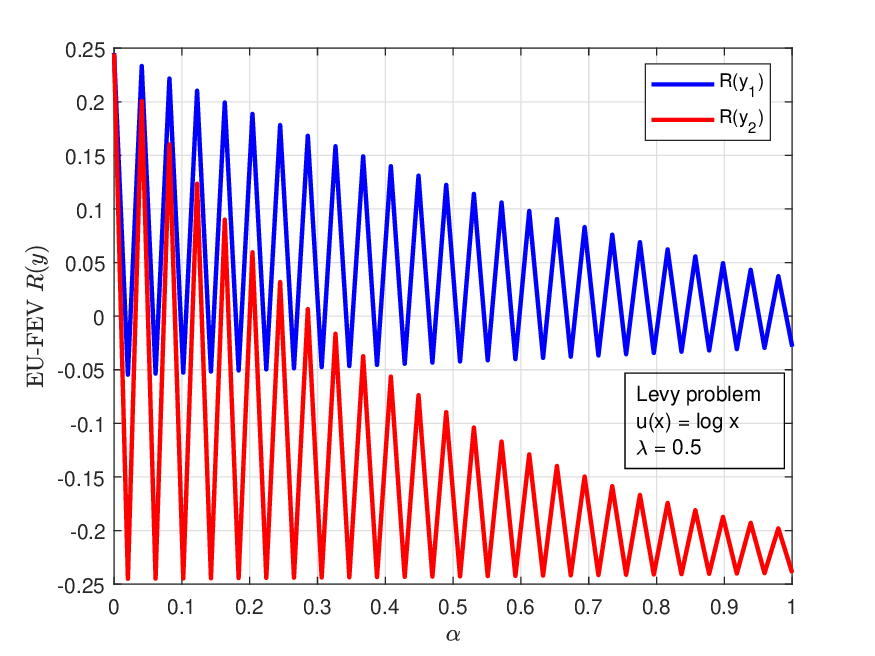} \\
\includegraphics[width=0.4\textwidth]{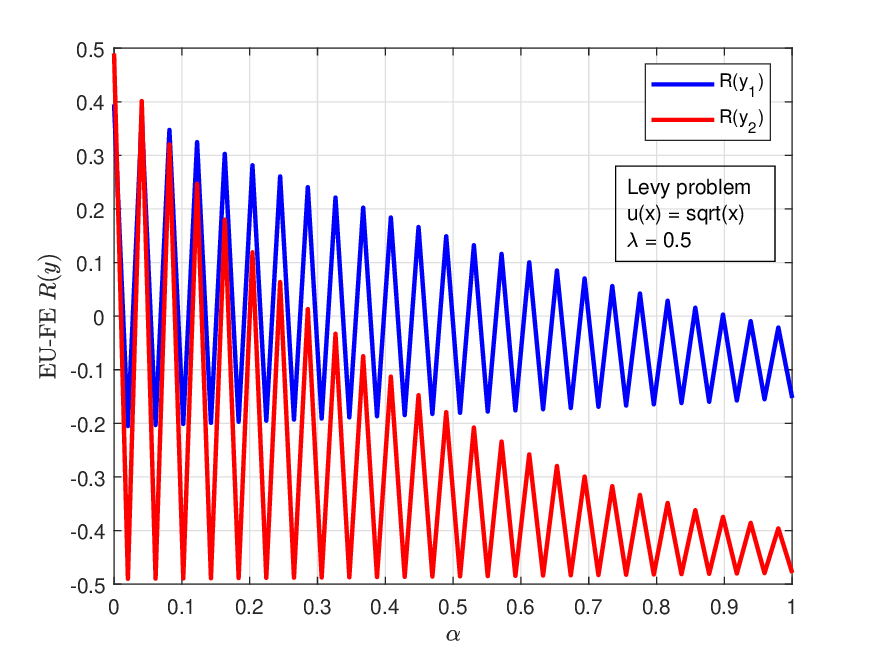} ~~ \includegraphics[width=0.4\textwidth]{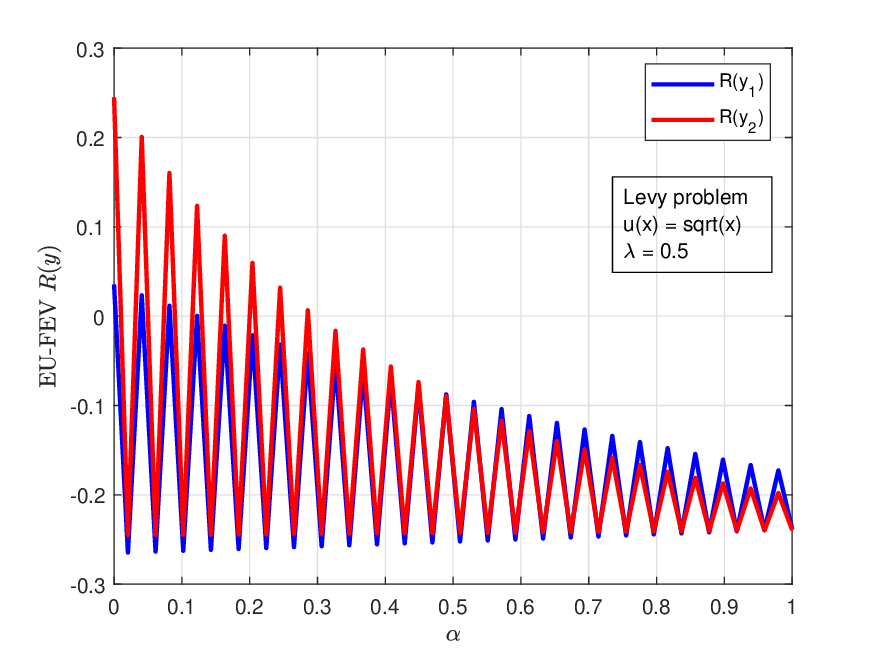}
\caption{Risk measures for actions $ y_3 $ and $ y_4 $ for different $ \alpha $ and $ u(x) $ in Levy problem}
\label{Fig10} 
\end{figure} 	
\begin{figure}[h]
\centering
\includegraphics[width=0.4\textwidth]{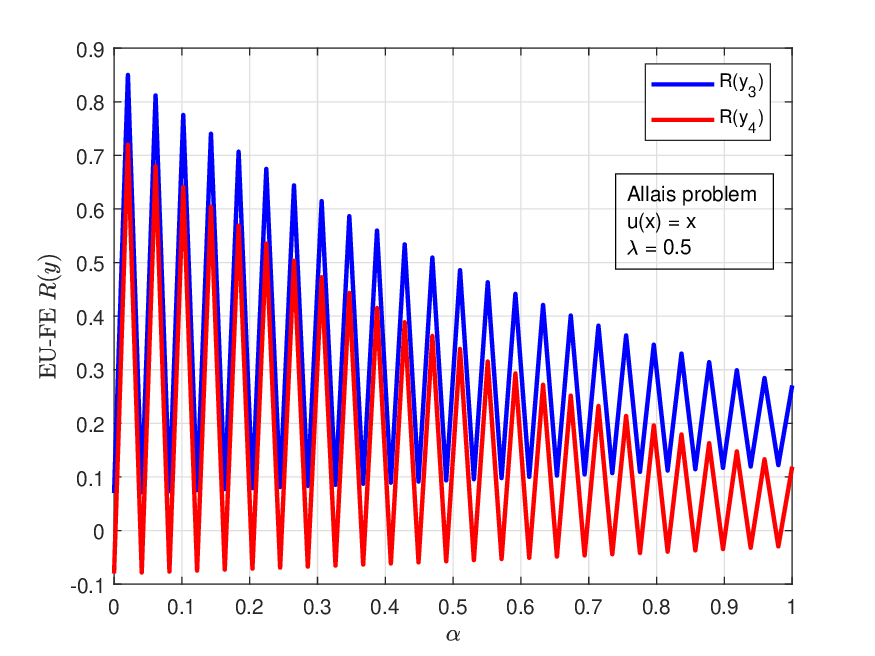} ~~ \includegraphics[width=0.4\textwidth]{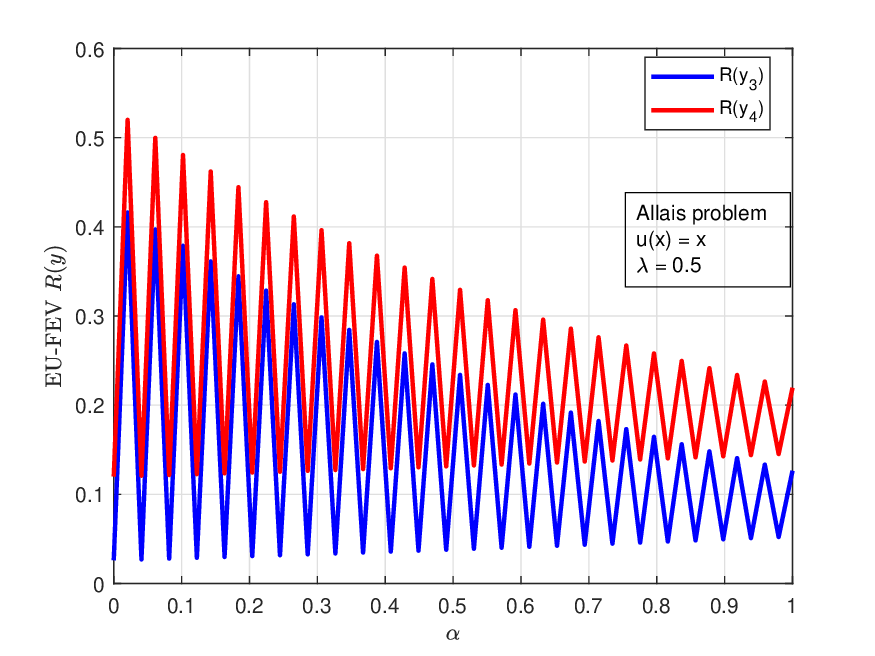} \\
\includegraphics[width=0.4\textwidth]{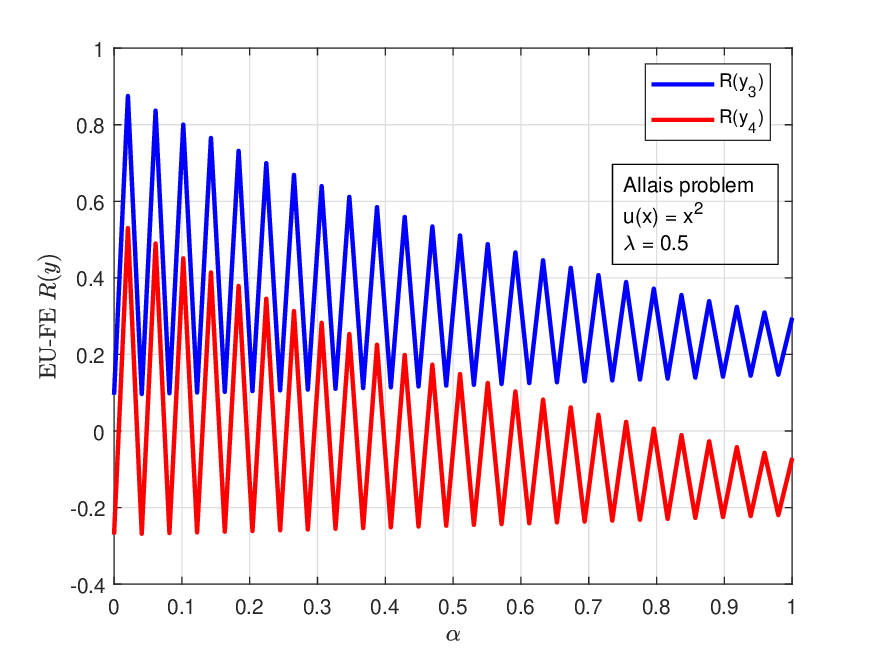} ~~ \includegraphics[width=0.4\textwidth]{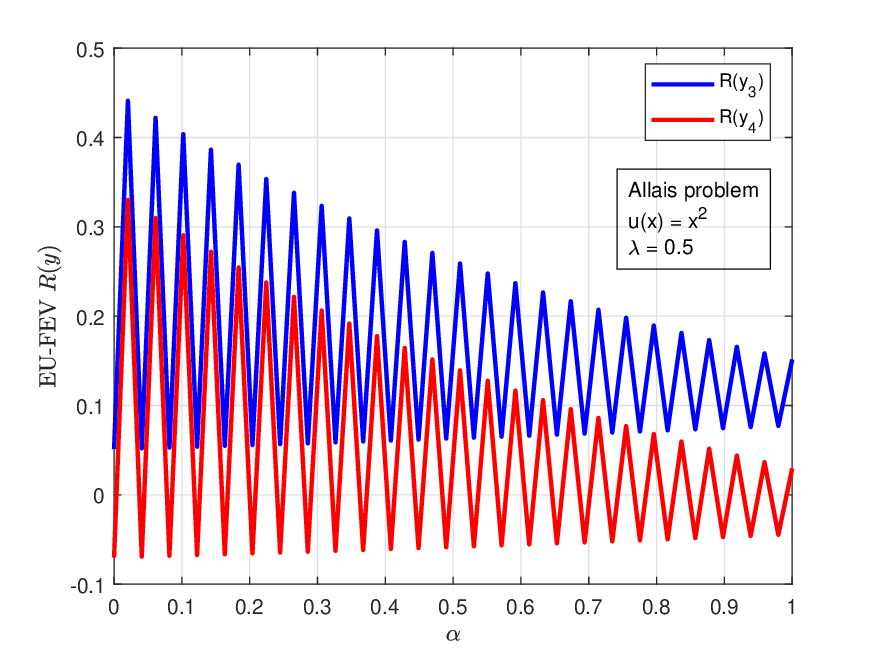} \\
\includegraphics[width=0.4\textwidth]{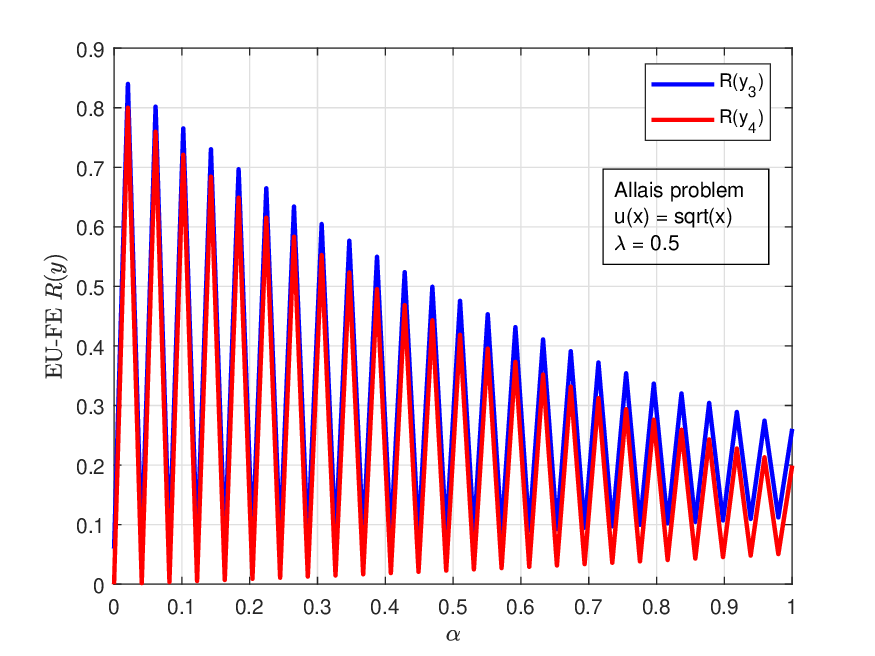} ~~ \includegraphics[width=0.4\textwidth]{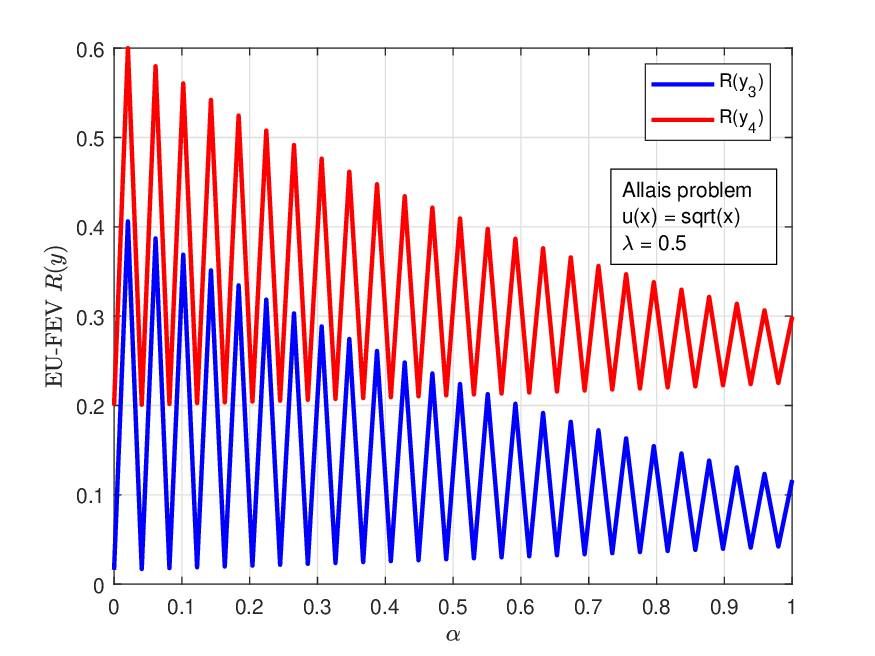}
\caption{Risk measures for actions $ y_3 $ and $ y_4 $ for different $ \alpha $ and $ u(x) $ in Allais problem}
\label{Fig11} 
\end{figure} 

\end{document}